\documentclass[a4paper]{article}
\usepackage{amssymb,amscd,amsfonts,mathrsfs,amsmath}
\usepackage[all]{xy}
\input amssym.def

\def\double{\mathbb}

\def\cc{{\double C}}
\def\nn{{\double N}}
\def\zz{{\double Z}}

\def\rr{{\double R}}

\newtheorem{theorem}{Theorem}[section]
\newtheorem{lemma}[theorem]{Lemma}
\newtheorem{corollary}[theorem]{Corollary}
\newtheorem{definition}[theorem]{Definition}
\newtheorem{proposition}[theorem]{Proposition}
\newtheorem{remark}[theorem]{Remark}

\def\res{\mathop{\mathrm{Res}}\limits_{z=0}}
\def\Pf{\mathop{\mathrm{Pf}}}

\def\cp{\rtimes}
\def\si{\sigma}

\def\cinf{C^{\infty}}
\def\cinfc{C^{\infty}_c}

\newcommand{\be}{\begin{equation}}
\newcommand{\ee}{\end{equation}}
\newcommand{\beq}{\begin{eqnarray}}
\newcommand{\eeq}{\end{eqnarray}}
\newcommand{\om}{\omega}
\newcommand{\Om}{\Omega}
\newcommand{\al}{\alpha}
\def\nat{\natural}
\def\id{{\mathop{\mathrm{id}}}}

\newcommand{\La}{\Lambda}
\newcommand{\la}{\lambda}
\newcommand{\Ec}{{\mathscr E}}

\newcommand{\non}{\nonumber}
\newcommand{\eps}{\varepsilon}

\newcommand{\Rc}{{\mathscr R}}

\newcommand{\Jc}{{\mathscr J}}

\newcommand{\Ind}{{\mathop{\mathrm{Ind}}}}

\def\re{\mathrm{Re}}

\newcommand{\Tr}{{\mathop{\mathrm{Tr}}}}

\newcommand{\tr}{{\mathop{\mathrm{tr}}}}
\newcommand{\Ac}{{\mathscr A}}

\newcommand{\cqfd}{\hfill\rule{1ex}{1ex}}

\def\Id{\mathrm{Id}}

\def\d{\partial}
\def\dd{\mathrm{\bf d}}

\def\Bc{{\mathscr B}}
\def\Cc{{\mathscr C}}
\def\Jc{{\mathscr J}}
\def\Kc{{\mathscr K}}
\def\Fc{{\mathscr F}}

\def\ker{\mathop{\mathrm{Ker}}}

\def\hom{{\mathop{\mathrm{Hom}}}}
\def\dom{{\mathop{\mathrm{Dom}}}}

\def\hotimes{\hat{\otimes}}

\def\Tt{\widetilde{T}}

\def\Omh{\widehat{\Omega}}
\def\Th{\widehat{T}}

\def\Rch{\widehat{\mathscr R}}

\def\Kc{{\mathscr K}}

\def\mod{\ \mathrm{mod}\ }

\def\supp{\mathrm{supp}\,}

\def\eh{\hat{e}}

\def\uh{\hat{u}}

\def\CL{\mathrm{CL}}
\def\LS{\mathrm{LS}}
\def\L{\mathrm{L}}
\def\CS{\mathrm{CS}}

\def\i{\mathrm{i}}

\def\R{\mathrm{R}}
\def\top{\mathrm{top}}

\def\Res{\mathrm{Res}}

\def\Ad{\mathrm{Ad}}

\begin{document}

\begin{center}

{\bf LOCAL INDEX THEORY FOR CERTAIN FOURIER INTEGRAL OPERATORS ON LIE GROUPOIDS}
\vskip 1cm
{\bf Denis PERROT}
\vskip 0.5cm
Institut Camille Jordan, CNRS UMR 5208\\
Universit\'e de Lyon, Universit\'e Lyon 1,\\
43, bd du 11 novembre 1918, 69622 Villeurbanne Cedex, France \\[2mm]
{\tt perrot@math.univ-lyon1.fr}\\[2mm]
\end{center}
\vskip 0.5cm
\begin{abstract}
We develop a local index theory for Fourier-integral operators associated to non-proper and non-isometric actions of Lie groupoids on smooth submersions. To such action is associated a short exact sequence of algebras, relating genuine Fourier-integral operators to their non-commutative symbol. We then compute the connecting map induced by this extension on periodic cyclic cohomology. When cyclic cohomology is localized at appropriate isotropic submanifolds of the groupoid in question, we find that the connecting map is expressed in terms of an explicit Wodzicki-type residue formula, which involves the jets of non-commutative symbols at the fixed-point set of the action.
\end{abstract}

\vskip 0.5cm

\noindent {\bf Keywords:} extensions, $K$-theory, cyclic cohomology.\\
\noindent {\bf MSC 2000:} 19D55, 19K56.

\section{Introduction}

This article concerns the index theory of a certain class of operators associated to smooth actions of Lie groupoids on manifolds. When a Lie groupoid $G$ acts \emph{properly} on a submersion of smooth manifolds $M\to B$, the index theory of a $G$-equivariant family of elliptic pseudodifferential operators on $M$ is rather well-understood \cite{CS, C94}. In his fundamental article \cite{C83}, Connes introduced cyclic cohomology techniques in order to deal with the $K$-theoretic index of a $G$-equivariant elliptic family. His result and subsequent generalizations use in a crucial manner the properness of the action, or in other circumstances, the fact that the action is isometric with respect to some Riemannian data on $M$. Much less is known about improper or non-isometric actions. Already in the much simpler case of a discrete group $G$ acting by diffeomorphisms on a closed manifold $M$ one wishes to study operators of the form
\be 
\sum_{g\in G} P_g U_g\ :\ \cinf(M)\to \cinf(M)\ ,\label{fou}
\ee
where $P_g$ is a pseudodifferential operator (say of order $\leq 0$) for any $g$, and $U_g$ is the representation of $g$ as a diffeomorphism on $M$. Such an operator is not pseudodifferential, and belongs to the larger class of Fourier-integral operators \cite{LNT}. Its principal symbol is not a smooth function on the cosphere bundle $S^*M$. Instead, it defines an element of the non-commutative crossed product algebra $\cinf(S^*M)\cp G$, or equivalently, the smooth convolution algebra of the \'etale groupoid $S^*M\cp G$. When $G$ is infinite this groupoid is not proper and the associated convolution algebra can be highly non-commutative. In the case of the group $G=\zz$, Savin and Sternin recently computed the index of an elliptic operator like (\ref{fou}) as a pairing between its leading symbol and a cyclic cohomology class of the algebra $\cinf(S^*M)\cp G$, see \cite{SS} and references therein. In the present article we want to unify these various index theorems and generalize them in two directions:
\begin{itemize}
\item Develop an index theory for operators whose non-commutative leading symbol belongs to the smooth convolution algebra of a \emph{not necessarily proper} groupoid such as the crossed product $S^*M\cp G$ above,
\item Evaluate the $K$-theoretical index of such operators on a wide range of cyclic cohomology classes, not necessarily localized at units.
\end{itemize}
The general situation considered in this article is the following. Let $G\rightrightarrows B$ be a Lie groupoid acting smoothly on a (surjective) submersion $\pi:M\to B$ of smooth manifolds. This basically means that any morphism $g\in G$ with source $s(g)\in B$ and range $r(g)\in B$ induces a diffeomorphism from the fiber $M_{r(g)}=\pi^{-1}(r(g))$ to the fiber $M_{s(g)}=\pi^{-1}(s(g))$, in a way compatible with the composition of morphisms in $G$. \emph{We do not impose any restriction on this action} (except smoothness), in particular, it is not necessarily proper nor isometric. Denote by $S^*_\pi M$ the bundle over $M$ whose fiber at a point $x\in M$ is the cotangent sphere of the submanifold $M_{\pi(x)}$ at $x$. Hence, $S^*_\pi M$ is the ``vertical" cosphere bundle over the submersion $M$. It is still endowed with a smooth action of $G$, and we consider the action groupoid $S^*_\pi M\cp G$. Its smooth convolution algebra 
\be 
\Ac=\cinfc(S^*_\pi M)\cp G
\ee 
is highly non-commutative in general. It may naturally be identified with the algebra of leading symbols of ``vertical non-commutative pseudodifferential operators'' on $M$. Let us explain what are these operators. Denote by $\CL^0_c(M)\to B$ the bundle with base $B$, whose fiber over a point $b\in B$ is the algebra of compactly supported \emph{classical} pseudodifferential operators of order $\leq 0$ on the manifold $M_b=\pi^{-1}(b)$. This bundle carries a natural action of $G$. Hence the algebra of smooth compactly supported sections of vertical classical pseudodifferential operators $\cinfc(B,\CL^0_c(M))$ can be twisted by the action of $G$. This leads to the crossed-product algebra of non-commutative pseudodifferential operators
\be 
\Ec=\cinfc(B,\CL_c^0(M))\cp G\ .
\ee
Clearly the projection of classical pseudodifferential operators of order $\leq 0$ onto the homogeneous component of degree 0 of their symbol extends to a surjective homomorphism of algebras $\Ec\to\Ac$, whose kernel is the two-sided ideal $\Bc=\cinfc(B,\CL_c^{-1}(M))\cp G$ of operators of order $\leq -1$ in $\Ec$. One thus gets a short exact sequence (extension) of algebras 
\be
(E): \ 0\to \Bc \to \Ec \to \Ac \to 0\ .
\ee
We say that an operator in the algebra $\Ec$ (unitalized) is \emph{elliptic} if its leading symbol is invertible in the algebra $\Ac$ (unitalized). This is a purely algebraic notion. An invertible leading symbol naturally defines an algebraic $K$-theory class in $K_1(\Ac)$. The \emph{index map} of the extension $(E)$ is the morphism
\be 
\Ind_E\ :\ K_1(\Ac)\to K_0(\Bc)
\ee
induced on algebraic $K$-theory in low degrees \cite{M}. Thus, if $[u]\in K_1(\Ac)$ is represented by the non-commutative symbol $u\in GL_{\infty}(\Ac)$ of an elliptic operator, its index $\Ind_E([u])$ is a $K$-theory element represented by an idempotent (matrix) in $\Bc$. Our goal is to evaluate this index on genuine cyclic cohomology classes of $\Bc$. In fact, as shown by Cuntz and Quillen in \cite{CQ97}, periodic cyclic cohomology satisfies excision in full generality. This means that the extension $(E)$ leads to a cohomology long exact sequence, with connecting map
\be 
E^*\ :\ HP^\bullet(\Bc) \to HP^{\bullet+1}(\Ac)\ .
\ee
Then Nistor \cite{Ni} (see also \cite{P8}) remarked that the index map $\Ind_E$ in algebraic $K$-theory is adjoint to the excision map $E^*$ with respect to the Chern-Connes pairing. Hence for any $[u]\in K_1(\Ac)$ and $[\tau]\in HP^0(\Bc)$ one has the equality of pairings
\be 
\langle [\tau], \Ind_E([u]) \rangle = \langle E^*([\tau]), [u] \rangle\ .
\ee
The left-hand side of this equality is a number, or ``higher index" caracterizing the $K$-theoretic index class of an elliptic operator. The right-hand side computes the higher index by means of a ``formula" involving the leading symbol of this operator. The difficulty is then to compute the image of cyclic cocycles $[\tau]$ under the excision map, which is actually not easy at all. We will use the general theory developed in \cite{P8, P9} for the explicit computation of the excision map in terms of \emph{local formulas}. The first step in this direction is to find which classes $[\tau]$ are ``good enough" to allow this computation. In our case, the algebra $\Bc = \cinfc(B,\CL^{-1}_c(M))\cp G$ is a finitely-summable thickening of the smooth convolution algebra of the pullback groupoid $\pi^*G\rightrightarrows M$. The latter is Morita equivalent to the groupoid $G\rightrightarrows B$, and the \emph{topological} cyclic cohomology of their respective convolution algebras $\cinfc(M)\cp \pi^*G$ and $\cinfc(B)\cp G$ are isomorphic. The topological cyclic cohomology of a smooth convolution algebra is defined according to its natural locally convex topology. We construct a canonical map (Proposition \ref{ptoto})
\be 
\tau_*\ :\ HP^\bullet_\top(\cinfc(B)\cp G) \to HP^\bullet(\cinfc(B,\CL^{-1}_c(M))\cp G)\ .  \label{tau}
\ee
By doing so we need to establish not so well-known properties of localization and Morita invariance of the periodic cyclic cohomology of smooth convolution algebras for general Lie groupoids. To our knowledge this was previously done only in the case of foliation groupoids (see \cite{BN, CrMo}), that is, groupoids which are Morita equivalent to \'etale ones. The following result shows that on the range of (\ref{tau}) the excision map $E^*$ factors through the topological cyclic cohomology of the convolution algebra $\Ac=\cinfc(S^*_\pi M)\cp G$:\\

\noindent{\bf Theorem 4.7} {\it Let $G\rightrightarrows B$ be a Lie groupoid and $\pi:M\to B$ a $G$-equivariant surjective submersion. Then one has a commutative diagram
\be 
\vcenter{\xymatrix{
HP^\bullet(\cinfc(B,\CL_c^{-1}(M))\cp G) \ar[r]^{\quad E^*} & HP^{\bullet+1}(\cinfc(S^*_\pi M)\rtimes G)  \\
HP_\top^{\bullet}(\cinfc(B)\cp G) \ar[u]^{\tau_*} \ar[r]^{\pi^!_G\quad} & HP_\top^{\bullet+1}(\cinfc(S^*_\pi M)\rtimes G) \ar[u] }}
\ee}

\noindent It remains to compute the map $\pi^!_G$. Let $O\subset G$ be an $\Ad$-invariant isotropic submanifold of $G$. This means that the range and source maps $O\rightrightarrows B$ coincide, and $O$ is globally invariant under the adjoint action of $G$. Then one has the notion of topological cyclic cohomology $HP_\top^{\bullet}(\cinfc(B)\cp G)_{[O]}$ localized at $O$, togetether with a forgetful map from the localized to the delocalized cohomology. Under suitable non-degeneracy hypotheses concerning the action of $O$ on $M$, we are able to calculate the map $\pi^!_G$ by means of explicit formulas involving \emph{residues of zeta-functions}. The use of zeta-function is motivated by the approach of Connes and Moscovici to the local index formula in non-commutative geometry \cite{CM95}. These residues generalize the well-known Wodzicki residue for classical pseudodifferential operators \cite{Wo}. They are given by integrals, over the cosphere bundle of the fixed point set for $O$, of certain local expessions in the complete symbol of the operators involved. As a refinement of Theorem \ref{ttop} we obtain the following result.\\ 

\noindent {\bf Theorem 5.6} {\it Let $G\rightrightarrows B$ be a Lie groupoid and let $O$ be an $\Ad$-invariant isotropic submanifold of $G$. Let $\pi:M\to B$ be a $G$-equivariant surjective submersion and assume the action of $O$ on $M$ non-degenerate. Then one has a commutative diagram 
\be 
\vcenter{\xymatrix{
HP^\bullet(\cinfc(B,\CL_{c}^{-1}(M))\cp G) \ar[r]^{\quad E^*} & HP^{\bullet+1}(\cinfc(S^*_\pi M)\rtimes G)  \\
HP_\top^{\bullet}(\cinfc(B)\cp G)_{[O]} \ar[u]^{\tau_*} \ar[r]^{\pi^!_G\qquad} & HP_\top^{\bullet+1}(\cinfc(S^*_\pi M)\rtimes G)_{[\pi^* O]} \ar[u] }}
\ee
where the isotropic submanifold $\pi^*O \subset S^*_\pi M\rtimes G$ is the pullback of $O$ by the submersion $S^*_\pi M\to B$. The map $\pi^!_G$ is given by an explicit residue formula.}\\

\noindent We will explain how to extract local index formulas from this theorem in a forthcoming paper. \\

Let us now give a brief description of the article. Section \ref{sgrou} recalls elementary notions about Lie groupoids, convolution algebras, and pseudodifferential operators. We introduce the basic extension associated to the action of a Lie groupoid on a submersion and describe the index map. In section \ref{scy} we review the cyclic cohomology of the smooth convolution algebra of a Lie groupoid, both in the algebraic and topological case, and establish the properties of localization and Morita invariance. In section \ref{sexci} we recall the computation of the excision map from \cite{P8} and prove Theorem \ref{ttop}. In section \ref{szeta} we use zeta-function renormalization techniques in order to prove the residue Theorem \ref{tres}.   \\

\noindent {\bf Acknowledgements:} This work was completed during a stay at the Institut de Math\'ematiques de Jussieu, Universit\'e Paris 7. The author is very grateful to all members of the team ``Alg\`ebres d'Op\'erateurs'' for their warm hospitality and the excellent working conditions, and especially to Georges Skandalis for many helpful discussions.

\section{Convolution algebras}\label{sgrou}

In this section we recall some basic facts about Lie groupoids, convolution algebras, pseudodifferential operators, and define the index of a non-commutative elliptic symbol as an algebraic $K$-theory class. 
\begin{definition}
A Lie groupoid $G\rightrightarrows B$ consists of:\\

\noindent a) Two smooth manifolds $B=G^{(0)}$ (the set of units) and $G=G^{(1)}$ (the set of morphisms). We will assume that $B$ and $G$ are Hausdorff and without boundary;\\
b) Two submersions $r,s:G\to B$ called the \emph{range} and \emph{source} map respectively;\\
c) A smooth map $m:G^{(2)}\to G$, where $G^{(2)}=\{(g_1,g_2)\in G\times G \ |\ s(g_1)=r(g_2)\}$ is the set of composable arrows, called the \emph{product} map. We usually write $m(g_1,g_2)= g_1g_2$;\\
d) A smooth embedding $u:B\hookrightarrow G$ and a diffeomorphism $i:G\to G$ called the \emph{unit} and \emph{inverse} map respectively. We usually write $u(b)=b$ for any $b\in B$ and $i(g)=g^{-1}$ for any $g\in G$.\\

\noindent These data are subject to compatibility conditions: for all composable morphisms $g_1, g_2, g_3 \in G$ one has\\

\noindent i) $r(g_1g_2)=r(g_1)$ and $s(g_1g_2)=s(g_2)$;\\
ii) $(g_1g_2)g_3 = g_1(g_2g_3)$ (associativity of the product);\\
iii) $g s(g)=g$ and $r(g)g=g$ (units);\\
iv) $g g^{-1}=r(g)$ and $g^{-1}g=s(g)$ (inverse).
\end{definition}
Following the standard convention we denote by $G^{(n)}$ the submanifold of composable $n$-tuples of morphisms $(g_1,\ldots,g_n)$ in $G^n$. One should keep in mind that a morphism $g\in G$ may be represented by a left-oriented arrow
$$
r(g) \stackrel{g}{\longleftarrow} s(g)
$$
and the product of morphisms is the concatenation of arrows. For any unit $b\in B$, we denote by $G_b=\{g\in G\ |\ s(g)=b\}$ the fiber of the source map over $b$, and by $G^b=\{g\in G\ |\ r(g)=b\}$ the fiber of the range map. These are submanifolds in $G$. At any point $b\in B$, the intersection $G_b^b = G_b\cap G^b$ is a Lie group called the isotropy group at $b$. The isotropy subset $I=\{g\in G\ |\ r(g)=s(g)\}$ is the union of all isotropy groups. It is a closed subset in $G$ but generally not a submanifold; the isotropy groups $G_b^b$ may have jumps when $b$ varies.  Let us recall some basic constructions. 

\begin{definition}
Let $G\rightrightarrows B$ be a Lie groupoid and $\pi: M\to B$ a submersion of smooth manifolds. The \emph{pullback groupoid} of $G$ by $\pi$ is the Lie groupoid
\be 
\pi^*G \rightrightarrows M
\ee
where $\pi^*G = \{(x,g,y)\in M\times G\times M\ |\ \pi(x)=r(g)\ ,\ \pi(y)=s(g)\}$, with range map $(x,g,y)\mapsto x$, source map $(x,g,y)\mapsto y$, and product of composable morphisms $(x,g_1,y)\cdot(y,g_2,z) = (x,g_1g_2,z)$.
\end{definition}

\begin{definition}
Two Lie groupoids $G_1\rightrightarrows B_1$ and $G_2\rightrightarrows B_2$ are \emph{Morita equivalent} if there exists a smooth manifold $M$ and two surjective submersions 
$$
B_1 \stackrel{\pi_1}{\longleftarrow} M  \stackrel{\pi_2}{\longrightarrow} B_2
$$
together with an isomorphism $\pi^*_1 G_1 \cong \pi^*_2 G_2$ between the pullback groupoids. Morita equivalence is an equivalence relation among Lie groupoids.
\end{definition}

\begin{definition}
Let $M$ be a smooth manifold and $G\rightrightarrows B$ a Lie groupoid. A right $G$-action on $M$ is given by \\
a) A smooth submersion $\pi:M\to B$;\\
b) A map from the fibered product $M\times_{(\pi,r)} G=\{(x,g)\in M\times G\ | \ \pi(x) = r(g)\}$ to $M$, sending a pair $(x,g)$ to the element $x\cdot g$ such that $\pi(x\cdot g)=s(g)$ and $(x\cdot g_1)\cdot g_2 =x\cdot(g_1g_2)$ whenever $g_1$ and $g_2$ can be composed.\\

\noindent We call such $\pi:M\to B$ a $G$-equivariant submersion.
\end{definition}
Note that the fibered product manifold $M\times_{(\pi,r)} G$ endowed with the composition law $(x,g_1)\cdot(y,g_2)=(x,g_1g_2)$ whenever $x\cdot g_1 = y$, is a Lie groupoid with range map $(x,g)\mapsto x$ and source map $(x,g)\mapsto x\cdot g$. We usually denote this \emph{action groupoid} by $M\cp G \rightrightarrows M$. It should not be confused with the pullback groupoid $\pi^* G$. \\
Let $M_b=\pi^{-1}(b)$ be the preimage of a point $b\in B$. Since $\pi$ is a submersion $M_b$ is a submanifold of $M$. Any element $g\in G$ induces a diffeomorphism $M_{r(g)} \to M_{s(g)}$ by $x\mapsto x\cdot g$. We extend this to a diffeomorphism $T^*M_{r(g)}\to T^*M_{s(g)}$ between the cotangent bundles of the respective submanifolds in $M$. One has a natural action of the multiplicative group $\rr^{\times}_+$ on the fibers of the cotangent bundle and the quotient $S^*M_{b}=T^*M_{b}/\rr^{\times}_+$ defines the \emph{cosphere bundle} over $M_{b}$. Moreover the diffeomorphism induced by $g$ commutes with the action of $\rr^{\times}_+$, hence descends to a diffeomorphism $S^*M_{r(\gamma)}\to S^*M_{s(\gamma)}$. The collection
\be
S^*_\pi M = \bigcup_{b\in B} S^*M_b
\ee
of vertical cosphere bundles is also clearly a $G$-equivariant submersion with base $B$.\\

We now define the smooth convolution algebra of a Lie groupoid $(r,s): G \rightrightarrows B$. The kernel $\ker s_*$ of the tangent map $s_*:TG\to TB$ is the vector bundle over $G$ whose vectors are tangent to the fibers $G_b=s^{-1}(b)$ of the source map. These vectors are the infinitesimal generators of the left multiplication of $G$ on itself. Since left and right multiplication commute, the fibers $(\ker s_*)_{g_1}$ and $(\ker s_*)_{g_2}$ are canonically isomorphic whenever $r(g_1)=r(g_2)$. The restriction of $\ker s_*$ to the submanifold of units $B\subset G$ yields a vector bundle $AG$ over $B$ called the \emph{Lie algebroid} of $G$, and one has a canonical identification
$$
r^*(AG) \cong \ker s_*
$$
of vector bundles over $G$. Hence, any section of $AG$ over $B$ can be pulled back by the rank map $r$ to a ``right-invariant" section of $\ker s_*$ over $G$. Passing to the dual bundle $A^*G$ of the Lie algebroid and taking the maximal exterior power, one thus gets an isomorphism between the line bundle $r^*(|\La^{\max} A^*G|)$ and the line bundle $|\La^{\max} (\ker s_*)^*|$ of 1-densities along the submanifolds $G_b = s^{-1}(b)$. The smooth convolution algebra of $G$ is then defined as the $\cc$-vector space
\be 
\cinfc(B)\cp G := \cinfc(G,r^*(|\La^{\max} A^*G|))
\ee
of smooth (complexified) sections of this line bundle. The product of two sections $a_1,a_2 \in \cinfc(B)\cp G$ evaluated on a point $g\in G$ is given by an integral over  all possible decompositions of $g$ into products $g_1g_2$,
\be
(a_1 a_2)(g) = \int_{g_1g_2=g} a_1(g_1)\, a_2(g_2)\ ,
\ee
where $a_1(g_1)a_2(g_2)\in |\La^{\max}A^*G|_{r(g_1)}\otimes |\La^{\max}A^*G|_{r(g_2)}$. The integral makes sense because $r(g_1)=r(g)$ is fixed, and when the point $g_2$ varies in $G_{s(g)}$ the line $|\La^{\max}A^*G|_{r(g_2)}$ runs over the fibers of the 1-density bundle over $G_{s(g)}$. This product is associative but not commutative in general. The convolution algebra is not unital unless $G$ is \emph{\'etale} (i.e. $r$ and $s$ are local diffeomorphisms, equivalently the Lie algebroid is reduced to its zero-section $B$) and $B$ is compact. In the latter case the unit $e$ of the convolution algebra is $e(g)=0$ for $g\notin B$ and $e(b)=1$ for $b\in B$. \\
By choosing a trivialisation of $|\La^{\max} A^*G|$, that is a nowhere vanishing section, one obtains by pullback a right-invariant section of the 1-density bundle, or equivalently \emph{smooth Haar system} on $G$. A choice of Haar system allows one to identify the convolution algebra of $G$ with the space of complex-valued functions with compact support on $G$,
$$
\cinfc(B)\cp G \cong \cinfc(G)
$$
and transfer the convolution product on the latter. Since a choice of Haar system is non-canonical, it is sometimes more convenient to use the completely canonical definition of the product given above.\\
Let $\R\to B$ is a $G$-equivariant associative algebra bundle, that is, a bundle of associative algebras over $B$ endowed with an action by isomorphisms $U_g:\R_{s(g)}\to \R_{r(g)}$, $\forall g\in G$, compatible with the products in $G$. We suppose that the sections of this bundle are endowed with some ``smooth" structure ensuring that the subsequent constructions make sense. Then we define the convolution algebra of $G$ twisted by the bundle $\R$ as the space of compactly supported sections
\be 
\cinfc(B,\R)\cp G := \cinfc(G,r^*(\R\otimes |\La^{\max} A^*G|))
\ee
endowed with a slight generalization of the above convolution product. For two sections $a_1,a_2 \in \cinfc(B,\R)\cp G$ we set
\be
(a_1 a_2)(g) = \int_{g_1g_2=g} a_1(g_1)\cdot U_{g_1} a_2(g_2)\ ,
\ee
where the isomorphism $U_{g_1}: \R_{s(g_1)}\otimes |\La^{\max} A^*G|_{s(g_1)} \to \R_{r(g_1)}\otimes |\La^{\max} A^*G|_{s(g_1)}$ acts on the fiber of $\R$ but not on the fiber of the density bundle. As an example let $\pi:M\to B$ be a $G$-equivariant submersion, and take $\R$ as the bundle whose fiber over $b\in B$ is the commutative algebra of smooth functions with compact support $\cinfc(M_b)$. Then 
$$
\cinfc(B,\R) \cp G \cong \cinfc(M)\cp G
$$
is canonically isomorphic to the smooth convolution algebra of the action groupoid $M\cp G$. We will sometimes use the notation $\cinf_p(B,\R)\cp G$ for the crossed product algebra of \emph{properly supported} sections of the bundle $r^*(\R\otimes |\La^{\max} A^*G|)$ over $G$.\\

Now let $\CL_c^m(M_b)$ be the space of \emph{compactly supported} classical (1-step polyhomogeneous) pseudodifferential operators of order $m\in \zz$, acting on the space of smooth functions with compact support on manifold $M_b$. Such a linear operator on $\cinfc(M_b)$ has distribution kernel with compact support in $M_b\times M_b$, and in a local coordinate system $(x,p)$ on $T^*M_b$ its symbol has an asymptotic expansion $\si(x,p)\sim\sum_{j\geq 0}\si_{m-j}(x,p)$, with $\si_{m-j}(x,p)$ a positively homogeneous function of degree $m-j\in\zz$ in the variable $p$. For any $g\in G$ and any $P\in \CL_c^m(M_{s(g)})$, the pushforward ${U}_gP U_g^{-1}\in \CL_c^m(M_{r(g)})$ is the adjoint action of the linear isomorphism $U_g: \cinf_c(M_{s(g)})\to \cinf_c(M_{r(g)})$ induced by the diffeomorphism $g$. We denote by $\CL_c(M)\to B$ the bundle over the base manifold $B$, whose fiber at a point $b\in B$ is the algebra $\CL_c(M_b)$. Hence $\CL_c(M)$ is a $G$-bundle. The subbundle $\L^{-\infty}_c(M)$ whose fiber is the algebra of smoothing operators with compact support, is a two-sided ideal in $\CL_c(M)$. The quotient $\CS_c(M_c)=\CL_c(M)/\L^{-\infty}_c(M)$ defines the algebra bundle of \emph{formal symbols} over $B$. We will essentially focus on the algebra bundle $\CL^0_c(M)$ of classical pseudodifferential operators of order $\leq 0$, and its two-sided ideal $\CL^{-1}_c(M)$. The quotient $\LS^0_c(M) = \CL^0(M)/\CL^{-1}(M)$ of \emph{leading symbols} is isomorphic to the bundle whose fiber over a point $b\in B$ is the commutative algebra $\cinfc(S^*M_b)$ of smooth compactly supported functions on the cosphere bundle of $M_b$. One thus has a commutative diagram of algebra bundles over $B$
\be
\vcenter{\xymatrix{
0 \ar [r] & \L^{-\infty}_{c}(M) \ar[r] \ar[d] & \CL^0_{c}(M) \ar[r] \ar@{=}[d] & \CS^0_{c}(M)  \ar[r] \ar[d] & 0  \\
0 \ar [r] & \CL^{-1}_{c}(M) \ar[r] & \CL^0_{c}(M) \ar[r] & \LS^0_c(M)  \ar[r] & 0 }} \label{ext}
\ee
where the rows are exact sequences. The left vertical arrow is an injection, whereas the right vertical arrow is a surjection. Since all bundles are $G$-bundles, one can form the corresponding convolution algebras by crossed product with the action of $G$. Notice there are canonical isomorphisms 
$$
\cinfc(B,\L^{-\infty}_{c}(M))\cp G \cong \cinfc(M)\cp \pi^*G\ , \quad \cinfc(B,\LS^0_c(M))\cp G \cong \cinfc(S^*_\pi M)\cp G
$$
where $\cinfc(M)\cp \pi^*G$ is the smooth convolution algebra of the pullback groupoid $\pi^*G\rightrightarrows M$, and $\cinfc(S^*_\pi M)\cp G$ is the smooth convolution algebra of the action groupoid $S^*_\pi M\cp G\rightrightarrows S^*_\pi M$. All arrows of (\ref{ext}) being $G$-equivariant homomorphisms of algebra bundles over $B$, one gets a commutative diagram
$$
\vcenter{\xymatrix{
0 \ar [r] & \cinfc(M) \cp \pi^*G \ar[r] \ar[d] & \cinfc(B,\CL^0_{c}(M))\cp G \ar[r] \ar@{=}[d] & \cinfc(B,\CS^0_{c}(M)) \cp G  \ar[r] \ar[d] & 0  \\
0 \ar [r] & \cinfc(B,\CL^{-1}_{c}(M)) \cp G \ar[r] & \cinfc(B,\CL^0_{c}(M)) \cp G \ar[r] & \cinfc(S^*_\pi M) \cp G  \ar[r] & 0 }}
$$
where the rows are short exact sequences of associative algebras. The bottom row is our main object of interest. We set $\Ac = \cinfc(S^*_\pi M)\cp G$,  $\Bc = \cinfc(B,\CL^{-1}_{c}(M)) \cp G$, $\Ec = \cinfc(B,\CL^0_{c}(M)) \cp G$ and consider the extension
\be 
(E)\ :\ 0\to \Bc\to\Ec\to\Ac\to 0\ .
\ee
The boundary map of Milnor \cite{M} associates to any algebraic $K$-theory class $[u]\in K_1(\Ac)$ an index $\Ind_E([u])\in K_0(\Bc)$. Recall that a class $[u]$ is represented by an element $u$ in the group $GL_{\infty}(\Ac)$ of invertible infinite matrices of the form $u=1+v$ with $v\in M_{\infty}(\Ac)$. We can choose two matrices $P$ and $Q$ with entries in the algebra $\Ec$ (with a unit adjoined), which project respectively to $u$ and its inverse $u^{-1}$. Then $P$ and $Q$ are inverse to each other modulo the ideal of matrices with entries in $\Bc$. The index $\Ind_E([u])\in K_0(\Bc)$ is represented by the difference of idempotents $[e]-[e_0]$ where
\be
e=\left(\begin{matrix}  1-(1-QP)^2 & Q(2-PQ)(1-PQ) \\ (1-PQ)P & (1-PQ)^2 \end{matrix}\right)\ ,\qquad e_0= \left(\begin{matrix}  1 & 0 \\ 0 & 0 \end{matrix}\right)\ .\label{idem}
\ee
\begin{definition}
Let $G\rightrightarrows B$ be a Lie groupoid and $\pi : M\to B$ a smooth $G$-equivariant submersion. Set $\Ac=\cinfc(S^*_\pi M)\cp G$ and $\Bc=\cinfc(B,\CL_c^{-1}(M))\cp G$. An invertible matrix $u\in GL_{\infty}(\Ac)$ is called an \emph{elliptic symbol}. Its index is the $K$-theory class
\be
\Ind_E([u])\in K_0(\Bc)\ ,
\ee
the image of $[u]\in K_1(\Ac)$ under the boundary map associated to the natural extension $(E):0\to \Bc\to\Ec\to\Ac\to 0$.
\end{definition}

\section{Cyclic homology}\label{scy}

We recall the basic notions of cyclic homology. Let $\Ac$ be an associative $\cc$-algebra. The space of noncommutative $n$-forms over $\Ac$ is $\Om^n\Ac = \Ac^+\otimes \Ac^{\otimes n}$ for all $n\geq 1$, where $\Ac^+=\Ac\oplus \cc$ is the algebra obtained by adjoining a unit. For $n=0$ one has $\Om^0\Ac = \Ac$. We write $a_0da_1\ldots da_n$ (reps. $da_1\ldots da_n$) for the generic element $a_0\otimes a_1 \otimes \ldots \otimes a_n$ (reps. $1 \otimes a_1 \otimes \ldots \otimes a_n$) in $\Om^n\Ac$. A differential $d:\Om^n\Ac \to \Om^{n+1}\Ac$ is defined by $d(a_0da_1\ldots da_n) = d(a_0da_1\ldots da_n)$ and $d(da_1\ldots da_n)=0$, and of course $d^2=0$. The direct sum $\Om\Ac= \bigoplus_{n\geq 0} \Om^n\Ac$ is gifted with the unique graded product satisfying the Leibniz rule $d(\om_1\om_2) = d\om_1\om_2 + (-1)^{n_1}\om_1 d\om_2$ for all $\om_i \in \Om^{n_i}\Ac$. This turns $\Om\Ac$ into a differential graded algebra. The Hochschild boundary operator $b:\Om^n\Ac\to\Om^{n-1}\Ac$ is defined by $b(\om da) = (-1)^{n-1}[\om,a]$ for all  $\om$ of degree $n-1$ and $a\in\Ac$. Equivalently
\beq
b(a_0da_1\ldots da_n)  &=& a_0a_1da_2\ldots da_n +  \sum_{i=1}^{n-1} (-1)^i a_0da_1\ldots d(a_ia_{i+1} )\ldots da_n \non\\
&& + (-1)^n a_na_0da_1\ldots da_{n-1}
\eeq
for all $a_0\in \Ac^+$ and $a_i\in\Ac$, $i\geq 1$. Let $\kappa = 1-(bd+db)$ be the Karoubi operator. One has $\kappa(\om da)=(-1)^nda\, \om$ for all $n$-form $\om$ and $a\in\Ac$. The Connes boundary operator $B:\Om^n \Ac\to \Om^{n+1}\Ac$ is defined by $B=(1+\kappa+\ldots+\kappa^n)d$, or equivalently
\be
B(a_0da_1\ldots da_n)  =  \sum_{i=0}^{n} (-1)^{ni} da_i\ldots da_nda_0 \ldots da_{i-1} \ .
\ee
One has $b^2=bB+Bb=B^2=0$, hence $\Om\Ac$ endowed with the operators $(b,B)$ is a bicomplex. By definition the cyclic homology $HC_\bullet(\Ac)$ is the homology of the following total complex with boundary $b+B$:
\be
\vcenter{\xymatrix{ \ar[d]_b & \ar[d]_b & \ar[d]_b   \\
\Om^2\Ac \ar[d]_b  &  \Om^1\Ac \ar[d]_b \ar[l]_B  & \Om^0\Ac \ar[l]_B  \\
\Om^1\Ac \ar[d]_b & \Om^0\Ac \ar[l]_B &  \\
\Om^0\Ac &  & }} \label{bic}
\ee
A cyclic homology class in $HC_n(\Ac)$ is therefore represented by an inhomogeneous differential form $\om_n+\om_{n-2}+\ldots \in \Om^n\Ac\oplus \Om^{n-2}\Ac\oplus \ldots$ which is closed in the sense $b\om_n+B\om_{n-2}=0$, $b\om_{n-2}+B\om_{n-4}=0$, etc... The obvious shift of degree two obtained by deleting the first column of the cyclic bicomplex gives rise to the periodicity operator $S: HC_n(\Ac)\to HC_{n-2}(\Ac)$. Now complete the space of differential forms $\Om\Ac$ by taking direct products instead of direct sums, and write $\Omh\Ac= \prod_{n\geq 0}\Om^n\Ac$. The operator $b+B$ extends to a well-defined boundary operator on the completed space. The \emph{periodic} cyclic homology of $\Ac$ is defined as the homology of this complex:
\be
HP_\bullet(\Ac) = H_\bullet(\Omh\Ac,b+B)\ .
\ee
By construction the periodic cyclic homology is $\zz_2$-graded, that is $HP_n(\Ac)\cong HP_{n+2}(\Ac)$ for all $n$. Since the complex $(\Omh\Ac,b+B)$ can be recovered as the projective limit (under the operation $S$) of the cyclic bicomplex of $\Ac$, the periodic and non-periodic cyclic homologies are related by a Milnor $\lim^1$ exact sequence
\be
0 \to {\varprojlim_S} ^1 HC_{\bullet-1}(\cinfc(G)) \to HP_\bullet(\cinfc(G)) \to \varprojlim_S HC_\bullet(\cinfc(G)) \to 0
\ee
The cyclic cohomology groups $HC^n(\Ac)$ are defined through the dual complex of (\ref{bic}) over $\cc$. One simply has to replace the vector space $\Om^n\Ac$ by its dual $\Om^n\Ac'=\hom(\Om\Ac,\cc)$ over $\cc$, and transpose the boundaries $(b,B)$. The periodicity operator $S:HC^n(\Ac)\to HC^{n+2}(\Ac)$ now raises the degree by two. Periodic cyclic cohomology is defined as the cohomology of the direct sum  $\Om\Ac'=\bigoplus_{n\geq 0}\Om^n\Ac'$, which is a $\zz_2$-graded complex once endowed with the transposed of the operator $b+B$:
\be
HP^\bullet(\Ac) = H^\bullet(\Om\Ac',b+B)\ .
\ee
The link between periodic and non-periodic cyclic cohomology is simpler than the case of homology, since $HP^\bullet(\Ac)$ is the inductive limit over $S$ of the groups $HC^\bullet(\Ac)$. There are obvious bilinear pairings $HC^n(\Ac)\times HC_n(\Ac)\to \cc$ and $HP^n(\Ac)\times HP_n(\Ac)\to \cc$.\\

Let $G\rightrightarrows B$ be a Lie groupoid. For notational convenience we suppose that a smooth Haar system on $G$ has been fixed, so the smooth convolution algebra $\Ac=\cinfc(B)\cp G$ is isomorphic to $\cinfc(G)$. The space of noncommutative $n$-forms $\Om^n\cinfc(G) = \cinfc(G)^+\otimes \cinfc(G)^{\otimes n}$ is a subspace of the smooth functions with compact support on the union of manifolds $G^{n+1}\cup G^n$. Indeed a generic $n$-form $a_0da_1\ldots da_n$, with $a_i\in\cinfc(G)$, is a function of $n+1$ points in $G$,
$$
(a_0da_1\ldots da_n)(g_0,g_1,\ldots,g_n) = a_0(g_0)a_1(g_1)\ldots a_n(g_n)\ ,
$$
while a $n$-form $da_1\ldots da_n$ is a function of $n$ points:
$$
(da_1\ldots da_n)(g_1,\ldots,g_n) = a_1(g_1)\ldots a_n(g_n)\ .
$$
Thus $\Om\cinfc(G)$ is actually a (complicated) subspace of the smooth compactly supported functions on the manifold $\bigcup_{n\geq 0}(G^{n+1}\cup G^n)$. Let $I=I^{(1)}$ be the isotropy subset of $G$, i.e. the set of morphisms $g\in G$ such that $r(g)=s(g)$. Following \cite{BN}, we define the set of loops $I^{(n)}\subset G^n$ for higher $n$ as
\be 
I^{(n)} = \{ (g_1,\ldots,g_n)\in G^n\ |\ s(g_i)=r(g_{i+1})\ \forall i<n\ ,\ s(g_n)=r(g_1)\}
\ee
We want to show that the periodic cyclic (co)homology of $\cinfc(G)$ can be localized, in the sense that the information of the cyclic bicomplex is entirely contained in the vicinity of the set of loops. To make it precise, let $\Om^n\cinfc(G)_0 \subset \Om^n\cinfc(G)$ be the subspace of functions vanishing on some neighborhood of the set of loops $I^{(n+1)}\cup I^{(n)} \subset G^{n+1}\cup G^n$. We define the quotient space of localized forms
\be 
\Om^n\cinfc(G)_{(I)} = \Om^n\cinfc(G)/\Om^n\cinfc(G)_0
\ee
Clearly the direct sum $\Om\cinfc(G)_0 = \bigoplus_{n\geq 0}\Om^n\cinfc(G)_0$ is stable by the boundary operators $b,B$ hence yields a subcomplex of the cyclic bicomplex. Consequently, the operators $b,B$ descend to the quotient $\Om\cinfc(G)_{(I)}$. This leads to a localized cyclic bicomplex. The localized periodic cyclic homology of $\cinfc(G)$ is defined accordingly, through the completion $\Omh\cinfc(G)_{(I)} = \prod_{n\geq 0} \Om^n\cinfc(G)_{(I)}$:
\be 
HP_\bullet (\cinfc(G))_{(I)} = H_\bullet (\Omh\cinfc(G)_{(I)}, b+B)
\ee
The localized periodic cyclic cohomology of $\cinfc(G)$ is defined by duality. We let $\Om^n\cinfc(G)_{(I)}'=\hom(\Om^n\cinfc(G)_{(I)},\cc)$ be the space of $\cc$-linear functionals on localized $n$-forms. This is exactly the set of linear maps $\varphi: \Om^n\cinfc(G)\to \cc$ vanishing on the subspace $\Om^n\cinfc(G)_0$. Then  the direct sum $\Om\cinfc(G)_{(I)}' =  \bigoplus_{n\geq 0} \Om^n\cinfc(G)_{(I)}'$ is a $\zz_2$-graded complex with boundary the transposed of $b+B$, whence
\be 
HP^\bullet (\cinfc(G))_{(I)} = H^\bullet (\hom(\Om\cinfc(G)_{(I)}, \cc), b+B)\ .
\ee

\begin{proposition}[Localization: algebraic case]\label{plocalg}
Let $G$ be any Lie groupoid with isotropy subset $I$. Then the projection of cyclic bicomplexes $\Om\cinfc(G) \to \Om\cinfc(G)_{(I)}$ induces isomorphisms
\be 
HP_\bullet (\cinfc(G)) \cong HP_\bullet (\cinfc(G))_{(I)} \ ,\quad HP^\bullet (\cinfc(G)) \cong HP^\bullet (\cinfc(G))_{(I)} \ . \label{localg}
\ee
\end{proposition}
{\it Proof:} Denote by $HC_\bullet(\cinfc(G))_0$ and $HP_\bullet(\cinfc(G))_0$ the cyclic homologies computed by the subcomplex $\Om\cinfc(G)_0$ of forms vanishing in the vicinity of loops. Our goal is to prove that $HP_\bullet(\cinfc(G))_0=0$. Then using the homology six-term exact sequence 
$$
\xymatrix{ HP_0(\cinfc(G))_0 \ar[r] & HP_0(\cinfc(G)) \ar[r] & HP_0(\cinfc(G))_{(I)} \ar[d] \\
HP_1(\cinfc(G))_{(I)} \ar[u] & HP_1(\cinfc(G)) \ar[l] & HP_1(\cinfc(G))_0 \ar[l] }
$$
associated to the exact sequence of complexes $0 \to \Omh\cinfc(G)_0 \to \Omh\cinfc(G) \to \Omh\cinfc(G)_{(I)} \to 0$ gives the first isomorphism in (\ref{localg}). We first show that the periodicity operator $S: HC_n(\cinfc(G))_0 \to HC_{n-2}(\cinfc(G))_0$ vanishes. Let $(U_i)_{i\in I}$ be a locally finite open covering of the space $B$ of units in $G$, where $I\subseteq \nn$ is an ordered, at most countable set (thus each compact subset of $B$ intersects a finite number of $U_i$'s). Let $(c_i)_{i\in I}$ be a partition of unity relative to this covering, in the sense that $c_i\in \cinfc(U_i)$ for all $i\in I$ and $\sum_{i\in I}c_i(x)^2=1$ for all $x\in B$. We view each function $c_i\in\cinfc(B)$ as a multiplier of the algebra $\cinfc(G)$: for any $a\in\cinfc(G)$ and $g\in G$ set
$$
(c_ia)(g) := c_i(r(g))\, a(g)\ ,\qquad (ac_i)(g) := a(g) \, c_i(s(g))\ .
$$
We use the partition of unity to build a map $\rho$ from $\cinfc(G)$ to the algebra of infinite matrices $M_{\infty}(\cinfc(G))$ as follows: for each $a\in\cinfc(G)$, the matrix element of $\rho(a)$ in position $(i,j)$ is $c_iac_j$. By the compactness of the support of $a$, only a finite number of matrix elements are non-zero. The condition $\sum_ic_i^2=1$ shows that $\rho$ is a homomorphism of algebras. Therefore, the induced map $\rho_*: \Om^n\cinfc(G)\to \Om^nM_{\infty}(\cinfc(G))$ composed with the trace of matrices yields a morphism of cyclic bicomplexes $\tr\rho_*:\Om\cinfc(G)\to \Om\cinfc(G)$. Explicitly for any $n$-form $a_0da_1\ldots da_n$ we have
$$
\tr\rho_*(a_0da_1\ldots da_n) = \sum_{i_0,\ldots,i_n} (c_{i_0}a_0c_{i_1})d(c_{i_1}a_1c_{i_2})\ldots d(c_{i_n}a_nc_{i_0})\ ,
$$
and similarly for $\tr\rho_*(da_1\ldots da_n)$. This morphism clearly restricts to a morphism of subcomplexes $\Om\cinfc(G)_0\to \Om\cinfc(G)_0$. Observe also that the $n$-form $(c_{i_0}a_0c_{i_1})d(c_{i_1}a_1c_{i_2})\ldots d(c_{i_n}a_nc_{i_0})$, viewed as a smooth function on $G^{n+1}$, has compact support consisting of multiplets $(g_0, g_1, \ldots ,g_n)$ such that $s(g_0)$ and $r(g_1)$ are in the support of $c_{i_1}$, $s(g_1)$ and $r(g_2)$ are in the support of $c_{i_2}$, and so on. Thus, if the supports of the $c_i$'s are small enough, the function $\tr\rho_*(a_0da_1\ldots da_n)$ can be localized to an arbitrary small neighborhood of the set of loops $I^{(n+1)}$ in $G^{n+1}$. In particular if $a_0da_1\ldots da_n$ belongs to the subspace $\Om^n\cinfc(G)_0$, one can always find a suitably fine covering of $B$ together with a partition of unity such that $\tr\rho_*(a_0da_1\ldots da_n)$ vanishes. The next step is to provide a homotopy between the homomorphism $\rho$ and the natural inclusion $\cinfc(G)\hookrightarrow M_{\infty}(\cinfc(G))$ in the upper left matrix position. Indeed, consider the isomorphism of algebras 
$$
M_{\infty}(\cinfc(G)) \cong \left( \begin{matrix} \cinfc(G) & \cc^{\infty}_{\textup{row}}\otimes \cinfc(G) \\ \cc^{\infty}_{\textup{col}}\otimes \cinfc(G) & M_{\infty}(\cinfc(G)) \end{matrix} \right)
$$
where $\cc^{\infty}_{\textup{row}}$ and $\cc^{\infty}_{\textup{col}}$ are respectively the spaces of infinite row and column matrices, with finitely many non-zero entries in $\cc$. According to this $2\times 2$ matrix notation one has two relevant homomorphisms $\rho^0,\rho^1:\cinfc(G)\to M_{\infty}(\cinfc(G))$ given by
$$
\rho^0(a)= \left(\begin{matrix} a & 0 \\ 0 & 0 \end{matrix} \right) \ ,\qquad \rho^1(a)= \left( \begin{matrix} 0 & 0 \\ 0 & \rho(a) \end{matrix} \right)\ .
$$
Then $\tr\rho^0_*$ is the identity map on $\Om\cinfc(G)$, while $\tr\rho^1_*=\tr\rho_*$. Let $u=(u_i)_{i\in I}$ be the infinite row with $u_i=c_i$, and $v=(v_i)_{i\in I}$ the infinite column with $v_i=c_i$. Note that $u$ and $v$ may have infinitely many non-zero entries. Nevertheless the scalar product $uv=\sum_i c_i^2 =1$ is well-defined, so that $vu$ is an idempotent matrix with infinitely many non-zero entries. From this one can produce a matrix $W$ and its inverse $W^{-1}$, which are both \emph{multipliers} of the algebra $ M_{\infty}(\cinfc(G))$: 
$$
W = \left(\begin{matrix} 0 & -u \\ v & 1-vu \end{matrix} \right)\ ,\qquad W^{-1} = \left(\begin{matrix} 0 & u \\ -v & 1-vu \end{matrix} \right)
$$
One has $\rho(a)=vau$ for all $a\in\cinfc(G)$, and the elements $\rho^0(a)$, $\rho^0(a)W$, $W^{-1}\rho^0(a)$ and $W^{-1}\rho^0(a)W=\rho^1(a)$ are all in $M_{\infty}(\cinfc(G))$. A classical argument using rotation matrices in $M_2(\cc)$ then shows that the homomorphisms $\rho^0$ and $\rho^1$ are homotopic after tensoring by $M_2(\cc)$. Hence the morphisms of cyclic bicomplexes $\tr\rho^0_*=\Id$ and $\tr\rho^1_*=\tr\rho_*$  induce the same maps in cyclic homology after stabilization by the periodicity operator $S$. Applying this to the subcomplex $\Om\cinfc(G)_0$ shows that $S:HC_n(\cinfc(G))_0 \to HC_{n-2}(\cinfc(G))_0$ coincides with $S\circ\tr\rho_*$. By virtue of the above observation, for any given $n$-cycle $\om$ there is a choice of open covering of $B$ with partition of unity so that $\tr\rho_*(\om)=0$. Hence $S=0$ on $HC_n(\cinfc(G))_0$ as claimed. Note that we cannot  apply this argument directly to the periodic cyclic homology $HP_\bullet(\cinfc(G))_0$, because a periodic cycle is an \emph{infinite} sequence of $n$-forms so we may find no suitable covering of $B$. Instead, we use the $\lim^1$ exact sequence
$$
0 \to {\varprojlim_S} ^1 HC_{\bullet-1}(\cinfc(G))_0 \to HP_\bullet(\cinfc(G))_0 \to \varprojlim_S HC_\bullet(\cinfc(G))_0 \to 0
$$ 
Since $S=0$, one has $\varprojlim^1  HC_{\bullet-1}(\cinfc(G))_0= 0$ and $\varprojlim  HC_{\bullet}(\cinfc(G))_0=0$, hence $ HP_{\bullet}(\cinfc(G))_0$ vanishes as required.\\
Passing to cohomology, we observe that, as a vector space over $\cc$, the \emph{non}-periodic cyclic cohomology $HC^n(\cinfc(G))_0$ is the algebraic dual of the space $HC_n(\cinfc(G))_0$. Hence, by transposition of the above result, the suspension operator $S:HC^n(\cinfc(G))_0 \to HC^{n+2}(\cinfc(G))_0$ vanishes as well as the inductive  limit $HP^\bullet(\cinfc(G))_0 = \varinjlim_S HC^\bullet(\cinfc(G))_0$. The second isomorphism in (\ref{localg}) then follows from the six-term exact sequence relating the periodic cyclic cohomology groups $HP^\bullet(\cinfc(G))_0$, $HP^\bullet(\cinfc(G))$ and $HP^\bullet(\cinfc(G))_{(I)}$. \hfill \cqfd\\

When an algebra $\Ac$ comes equipped with a locally convex topology, the algebraic cyclic (co)homologies $HP_\bullet(\Ac)$ and $HP^\bullet(\Ac)$ described above can be replaced by appropriate topological versions. The topological cyclic homology of such an algebra is defined through a space of noncommutative differential forms as in the algebraic case, the only difference is that one has to replace algebraic tensor products by topological ones. The space of topological $n$-forms is thus $\Om_\top^n\Ac = \Ac^+\hotimes \Ac^{\hotimes n}$ where $\Ac^+=\Ac\oplus\cc$ is the algebra obtained by adjoining a unit, and $\hotimes$ is an appropriate completion of the algebraic tensor product. In the case of a Lie groupoid $G$, its smooth convolution algebra $\Ac=\cinfc(G)$ has the topology of an LF-space, which is the inductive limit topology, over all compact subsets $K\subset G$, of the Fr\'echet spaces $\cinf_K(G)$ of smooth functions with support contained in $K$.   In this case one may choose Grothendieck's inductive tensor product. Recall that if $M$ and $N$ are two manifolds, the inductive tensor product of LF-spaces $\cinfc(M)\hotimes\cinfc(N)$ is isomorphic to $\cinfc(M\times N)$. The space of noncommutative $n$-forms over the convolution algebra is thus isomorphic to
\be
\Om_\top^n\cinfc(G) \cong \cinfc(G^{n+1} \cup G^n)\ .
\ee
The operators $(b,B)$ extend by continuity to well-defined boundary operators on the direct sum $\Om_\top\cinfc(G) = \bigoplus_{n\geq 0} \Om_\top^n\cinfc(G)$ and also on the direct product $\Omh_\top\cinfc(G) = \prod_{n\geq 0} \Om_\top^n\cinfc(G)$. By definition the topological periodic cyclic homology of the convolution algebra is
\be
HP_\bullet^\top(\cinfc(G)) = H_\bullet (\Omh_\top\cinfc(G), b+B)\ .
\ee
Passing to cohomology one has to consider an appropriate dual space to noncommutative forms $\Om^n_\top\cinfc(G)'=\hom(\Om^n_\top\cinfc(G),\cc)$. We take the space of \emph{continuous}  and linear functionals $\varphi:\Om^n_\top\cinfc(G)\to\cc$ with \emph{bounded singularity order}. Such a functional $\varphi$ is exactly represented by a distribution on the manifold $G^{n+1}\cup G^n$ whose singularity order is \emph{finite}, say $k$: its evaluation on a smooth function $\om\in \Om^n_\top\cinfc(G)$ formally reads
$$
\varphi(\om) = \int_{G^{n+1}} {\varphi}_{0,n}(g_0,\ldots,g_n)\om(g_0,\ldots,g_n) + \int_{G^n} {\varphi}_{1,n}(g_1,\ldots,g_n)\om(g_1,\ldots,g_n)
$$
and can be extended to a continuous linear functional on functions of class $C^k$. We endow the space $\Om^n_\top\cinfc(G)'$ with the weak-$*$ topology. The transposed of the total operator $b+B$ acting on the direct sum $\Om_\top\cinfc(G)'=\bigoplus_{n\geq 0} \Om^n_\top\cinfc(G)'$ thus yields a $\zz_2$-graded topological complex. The topological periodic cyclic cohomology of the convolution algebra is by definition
\be
HP^\bullet_\top(\cinfc(G)) = H^\bullet (\Om_\top\cinfc(G)', b+B)\ .
\ee
We want to discuss localization of topological periodic cyclic (co)homology. Thus let $\Om^n_\top\cinfc(G)_0 \subset \Om^n_\top\cinfc(G)$ be the subspace of functions vanishing on some neighborhood of the set of loops $I^{(n+1)}\cup I^{(n)}$. This subspace is not closed in $\Om^n_\top\cinfc(G)$. We define the localized noncommutative $n$-forms as the (non-Hausdorff) quotient space
\be
\Om^n_\top\cinfc(G)_{(I)} = \Om^n_\top\cinfc(G)/\Om^n_\top\cinfc(G)_0\ .
\ee
Since the operators $(b,B)$ descend on localized forms, we define the localized topological periodic cyclic homology of the convolution algebra as the homology of the direct product $\Omh_\top\cinfc(G)_{(I)} = \prod_{n\geq 0} \Om_\top^n\cinfc(G)_{(I)}$
\be
HP_\bullet^\top(\cinfc(G))_{(I)} = H_\bullet (\Omh_\top\cinfc(G)_{(I)}, b+B)\ .
\ee
The definition of localized topological cyclic cohomology is analogous to algebraic setting. Let $\Om^n_\top\cinfc(G)_{(I)}'$ be the set of continuous bounded linear functionals $\varphi:\Om^n_\top\cinfc(G)\to\cc$ vanishing on the subspace $\Om^n_\top\cinfc(G)_0$ (hence also on its closure). These functionals are characterized by their distribution kernel whose support is entirely contained in the set of loops $I^{(n+1)}\cup I^{(n)}$. The direct sum $\Om_\top\cinfc(G)_{(I)}'=\bigoplus_{n\geq 0} \Om^n_\top\cinfc(G)_{(I)}'$ endowed with the transposed of $b+B$ is therefore a $\zz_2$-graded subcomplex of $\Om_\top\cinfc(G)'$ and we set
\be
HP^\bullet_\top(\cinfc(G))_{(I)} = H^\bullet (\Om_\top\cinfc(G)_{(I)}', b+B)\ .
\ee

\begin{proposition}[Localization: topological case]\label{ploctop}
Let $G$ be a Lie groupoid. The projection of cyclic bicomplexes $\Om_\top\cinfc(G) \to \Om_\top\cinfc(G)_{(I)}$ induces an isomorphism in periodic cyclic homology
\be 
HP^\top_\bullet (\cinfc(G)) \cong HP^\top_\bullet (\cinfc(G))_{(I)} \ . \label{loctop}
\ee
Moreover if the closure of $\Om_\top\cinfc(G)_0$ is a direct summand in $\Om_\top\cinfc(G)$ as a topological vector subspace, then the natural map in periodic cyclic cohomology
\be 
HP_\top^\bullet (\cinfc(G))_{(I)} \to HP_\top^\bullet (\cinfc(G))
\ee
is surjective (with kernel a topological vector space which does not separate zero from any other vector).
\end{proposition}
{\it Proof:} For cyclic homology, the proof of Proposition \ref{plocalg} applies verbatim. A partition of unity $(c_i)_{i\in I}$ relative to an open covering of $B=G^{(0)}$ yields a continuous homomorphism $\rho:\cinfc(G)\to M_{\infty}(\cinfc(G))$ by setting $\rho(a)_{ij} = c_iac_j$. The resulting chain map $\tr\rho_*: \Om_\top\cinfc(G)\to \Om_\top\cinfc(G)$ is explicitly described as follows. Any $n$-form $\om\in\Om^n_\top\cinfc(G)$ may be viewed as a smooth function over $G^{n+1}\cup G^n$. One has
\beq
\lefteqn{\big(\tr\rho_*(\om)\big)(g_0,\ldots ,g_n) = \om(g_0,\ldots,g_n) \times} \non\\ 
&& \sum_{i_0,\ldots,i_n} c_{i_0}(r(g_0))c_{i_1}(s(g_0))c_{i_1}(r(g_1))c_{i_2}(s(g_1))\ldots c_{i_n}(r(g_n))c_{i_0}(s(g_n))\ ,\non
\eeq
and similarly on $(g_1,\ldots,g_n)$. This expression vanishes if $\om$ belongs to the subspace $\Om^n_\top\cinfc(G)_0$ and the covering of $B$ is fine enough. The isomorphism $HP^\top_\bullet (\cinfc(G)) \cong HP^\top_\bullet (\cinfc(G))_{(I)}$ thus follows from homotopy invariance as before.\\
The case of cohomology requires some care, because we can no longer use a duality argument as in the algebraic setting. If we assume that the closure of $\Om_\top\cinfc(G)_0$ is a topological direct summand in $\Om_\top\cinfc(G)$, then at the dual level $\Om_\top\cinfc(G)'$ endowed with the weak-$*$ topology splits as the direct sum of $\Om_\top\cinfc(G)_{(I)}'$ and a closed supplementary subspace. In particular the short exact sequence of cyclic bicomplexes $0 \to \Om_\top\cinfc(G)_{(I)}' \to \Om_\top\cinfc(G)' \to \Om_\top\cinfc(G)_0'\to 0$ has a \emph{continuous} linear section. This implies the existence of a six-term exact sequence with continuous maps in topological periodic cyclic cohomology:
$$
\xymatrix{ HP_\top^0(\cinfc(G))_0 \ar[d] & HP_\top^0(\cinfc(G)) \ar[l] & HP_\top^0(\cinfc(G))_{(I)} \ar[l] \\
HP_\top^1(\cinfc(G))_{(I)} \ar[r] & HP_\top^1(\cinfc(G)) \ar[r] & HP_\top^1(\cinfc(G))_0 \ar[u] }
$$
The quotient complex $\Om_\top\cinfc(G)_0'$ is the continuous dual of the \emph{closure} of $\Om_\top\cinfc(G)_0$ inside $\Om_\top\cinfc(G)$. We will show that the periodic cyclic cohomology $HP_\top^\bullet(\cinfc(G))_0$ is degenerate in the sense that its topology does not separate zero from any other element. For this we endow the manifold $B$ with a riemannian metric. Consider a sequence of real numbers $\eps>0$ with limit $\eps\to 0$. For each $\eps$, we can choose an open covering $(U^\eps_i)_{i\in I}$ of $B$ together with a partition of unity $(c^\eps_i)_{i\in I}$ with the following properties: over any compact subset $K\subset B$ hold\\
i) For all $x\in K$, the number of $U_i$'s containing $x$ is bounded uniformly with respect to $x$ and $\eps$;\\
ii) The radius of $U^\eps_i\cap K$ is $\leq \eps$ for all $i\in I$ and $\eps$;\\
iii) The partial derivatives $\d_x^\al c^\eps_i$ are bounded by $C_{K,\al}\eps^{-|\al|}$ for all $i\in I$, $\eps>0$ and multi-index $\al=(\al_1,\ldots,\al_n)$. Here $n=\dim B$, $|\al|=\al_1+\ldots+\al_n$ and $C_{K,\al}$ is a constant independent of $\eps$.\\
These data give rise to a homomorphism $\rho^\eps:\cinfc(G)\to M_{\infty}(\cinfc(G))$ for each $\eps$ of the sequence, together with the associated chain map $\tr\rho^\eps_*:\Om_\top\cinfc(G)\to \Om_\top\cinfc(G)$. If a $n$-form $\om$ is in the closure of $\Om^n_\top\cinfc(G)_0$, then as a smooth function on $G^{n+1}\cup G^n$, $\om$ vanishes as well as all its partial derivatives on the set of loops $I^{(n+1)}\cup I^{(n)}$. This means that at a distance $\eps$ away from the set of loops, the partial derivatives of $\om$ grow slower than any power of $\eps$. From points i), ii), iii) above it follows that $\lim_{\eps\to 0} \tr\rho^\eps_*(\om)=0$. Passing to the continuous dual, any periodic cyclic cocycle $\varphi\in \Om_\top\cinfc(G)_0'$ is cohomologous to the sequence of cocycles $\varphi\circ\tr\rho^\eps_*$ which tends to zero in the weak-$*$ topology. Hence $HP_\top^\bullet(\cinfc(G))_0$ is a degenerate topological vector space as claimed. Now let $V$ be the range of the map $p:HP_\top^\bullet(\cinfc(G)) \to HP_\top^\bullet(\cinfc(G))_0$. The six-term exact sequence of periodic cyclic cohomology gives rise to a short exact sequence
$$
HP_\top^\bullet(\cinfc(G))_{(I)} \to HP_\top^\bullet(\cinfc(G)) \stackrel{p}{\to} V \to 0
$$
The topology of $V$ as a closed vector subspace of $HP_\top^\bullet(\cinfc(G))_0$ coincides with the quotient topology of $HP_\top^\bullet(\cinfc(G))/\ker p$. We know that it is degenerate, hence $\ker p$ must be dense in $HP_\top^\bullet(\cinfc(G))$. Since $p$ is continuous $\ker p$ is also necessarily closed. Hence one has $\ker p=HP_\top^\bullet(\cinfc(G))$, and the six-term exact sequence reduces to $0\to HP_\top^\bullet(\cinfc(G))_0 \to HP_\top^\bullet(\cinfc(G))_{(I)} \to HP_\top^\bullet(\cinfc(G)) \to 0$. \hfill\cqfd\\

\begin{remark}\label{rloc}
\textup{The proof of \ref{ploctop} shows that any topological cyclic cohomology class $[\varphi]\in HP^\bullet_\top(\cinfc(G))$ can be represented by a finite collection of distributions $\varphi:\Om^n_\top\cinfc(G)\to\cc$ whose supports are contained in an arbitrarily small open neighborhood of the localization set $I^{(n+1)}\cup I^{(n)}$. In general we do not know whether $[\varphi]$ can be represented by distributions with support exactly contained in $I^{(n+1)}\cup I^{(n)}$, unless the closure of $\Om_\top\cinfc(G)_0$ admits a topological supplementary subspace in $\Om_\top\cinfc(G)$. A sufficient condition for this to be true is that the space of loops $I^{(n)}$ is a smooth submanifold of $G^n$ for all $n$. This happens if the foliation $(B,\Fc)$ induced on the unit space of the groupoid $G\rightrightarrows B$ is non-singular, and is always verified, for example, by \'etale groupoids. This condition can be slightly relaxed by requiring that $I^{(n)}$ is a union of smooth submanifolds, with ``sufficiently nice" crossings. }
\end{remark}
The fact that the quotient $\Om^n_\top\cinfc(G)_{(I)}$ is a non-Hausdorff space is inconvenient. In order to deal with a much nicer space we introduce a strict localization of differential forms, quotienting by the closure of $\Om^n_\top\cinfc(G)_0$ inside $\Om^n_\top\cinfc(G)$:
\be
\Om^n_\top\cinfc(G)_{[I]} = \Om^n_\top\cinfc(G)/\, \overline{\Om^n_\top\cinfc(G)_0}\ .
\ee
We define accordingly the localized periodic cyclic cohomology $HP^\bullet_\top(\cinfc(G))_{[I]}$. Since the distributions vanishing on the subspace $\Om^n_\top\cinfc(G)_0$ and its closure actually coincide, one always has an isomorphism
\be 
HP^\bullet_\top(\cinfc(G))_{[I]} \cong HP^\bullet_\top(\cinfc(G))_{(I)}\ .
\ee
More generally if an isotropic subset $O\subset I$ is invariant under the adjoint action of $G$, we define $O^{(n)}\subset G^n$ as the set of composable arrows $(g_1,\ldots,g_n)$ with product $g_1\ldots g_n\in O$, and $\Om^n_\top\cinfc(G)_{0}^O$ as the functions with compact support on $G^{n+1}\cup G^n$ vanishing in a neighborhood of $O^{(n+1)}\cup O^{(n)}$. Then the strict localization at $O$ 
\be 
\Om^n_\top\cinfc(G)_{[O]} = \Om^n_\top\cinfc(G)/\,\overline{\Om^n_\top\cinfc(G)_{0}^O}
\ee
is a quotient complex of the cycic bicomplex. The corresponding localized cyclic cohomology $HP^\bullet_\top(\cinfc(G))_{[O]}$ is the cohomology of the complex of distributions with bounded singularity order, whose support is contained in $O^{(n+1)}\cup O^{(n)}$. Note that if $O^{(n)}$ is a submanifold of $G^n$ for all $n$, the quotient $\Om^n_\top\cinfc(G)_{[O]}$ is the space of jets of functions to any order at the localization manifold, i.e. the space of Taylor expansions of functions in the direction transverse to $O^{(n+1)}\cup O^{(n)}$. \\

We now study the invariance of cyclic homology with respect to Morita equivalences. The following lemma is based on an idea of G. Skandalis.
\begin{lemma}\label{lgeorges}
Let $G\rightrightarrows B$ be a Lie groupoid. Let $V\subset B$ be an open subset which intersects each orbit of $G$ and denote by $G_V\rightrightarrows V$ the restriction groupoid. Then one has isomorphisms in algebraic and topological periodic cyclic (co)homology
\beq
HP_\bullet (\cinfc(G_V)) &\cong& HP_\bullet (\cinfc(G)) \ ,\quad HP^\bullet (\cinfc(G_V)) \cong HP^\bullet (\cinfc(G)) \non\\
HP^\top_\bullet (\cinfc(G_V)) &\cong& HP^\top_\bullet (\cinfc(G)) \ ,\quad HP_\top^\bullet (\cinfc(G_V)) \cong HP_\top^\bullet (\cinfc(G)) \non
\eeq
induced by the inclusion of convolution algebras $\cinfc(G_V)\hookrightarrow \cinfc(G)$.
\end{lemma}
{\it Proof:} We first observe that the group of bisections of $G$ acts by multipliers on the convolution algebra $\cinfc(G)$. Indeed if $\beta:B\to G$ is a bisection, its left and right actions on an element $a\in \cinfc(G)$ are defined by
$$
(\beta\cdot a)(g) = a(\beta^{-1}(r(g))\cdot g)\ ,\qquad (a\cdot \beta)(g) = a(g\cdot \beta(s(g))^{-1})
$$
for all $g\in G$. One checks that the usual relations $\beta_1\cdot(\beta_2\cdot a)=\beta_1\beta_2\cdot a$, $(\beta\cdot a_1)a_2=\beta\cdot a_1a_2$ etc... are fulfilled, i.e. the group of bisections acts by multipliers on $\cinfc(G)$. If $\beta:U\to G$ is only a \emph{local} bisection over an open subset $U\subset B$, the actions $\beta\cdot a$ and $a\cdot\beta$ are still defined provided that the support of $a$ satisfies appropriate compatibility conditions with respect to the domain and range of the local diffeomorphism $\phi_\beta$ associated to $\beta$.\\
Now let $(\beta_i,U_i)_{i\in I}$ be a collection of local bisections $\beta_i:U_i\to G$ indexed by an at most countable set $I$, such that: i) the collection $(U_i)_{i\in I}$ is a locally finite covering of $B$ and ii) the range of the local diffeomorphism $\phi_{\beta_i}:U_i\to V_i$ is contained in the open subset $V\subset B$ for all $i\in I$. Condition ii) can be satisfied because $V$ intersects each orbit of $G$ by hypothesis. Then choose a partition of unity $(c_i)_{i\in I}$, with $\sum_ic_i(x)^2=1$, relative to the covering $(U_i)$. One builds an algebra homomorphism 
$$
\rho: \cinfc(G) \to M_{\infty}(\cinfc(G_V))
$$
by setting the $(i,j)$ entry of the matrix $\rho(a)$ equal to $\rho(a)_{ij}=\beta_i c_i a c_j \beta_j^{-1}$ for all $a\in\cinfc(G)$. Here the functions $c_i$ are multipliers of the convolution algebra as in the proof of Proposition \ref{plocalg}. An explicit computation gives, by evaluating on a point $g\in G$,
$$
\rho(a)_{ij}(g) = c_i\big(\phi_{\beta_i^{-1}}(r(g))\big)\, a\big(\beta_i^{-1}(r(g))\cdot g\cdot \beta_j^{-1}(s(g))^{-1}\big) \, c_j\big(\phi_{\beta_j^{-1}}(s(g))\big)\ .
$$
Since $\supp(c_i)\subset U_i$ and $\supp(c_j)\subset U_j$, the latter expression vanishes unless $\phi_{\beta_i^{-1}}(r(g))\in U_i$ and $\phi_{\beta_j^{-1}}(s(g))\in U_j$, that is, unless $r(g)\in V_i$ and $s(g)\in V_j$. This shows that $\rho(a)_{ij}$ is indeed an element of the subalgebra $\cinfc(G_V)$. The map induced by the homomorphism $\rho$ on differential forms composed with the trace map yields a morphism of cyclic bicomplexes $\tr\rho_*:\Om\cinfc(G)\to\Om\cinfc(G_V)$. We want to show that the latter is an isomorphism in periodic cyclic homology. The obvious candidate for an inverse comes from the morphism of cyclic bicomplexes $\iota_*:\Om\cinfc(G_V)\to\Om\cinfc(G)$ induced by the inclusion homomorphism $\iota:\cinfc(G_V)\to\cinfc(G)$. Hence it remains to prove that $\iota_*\circ \tr\rho_*$ and $\tr\rho_*\circ\iota_*$ are chain homotopic to the identity maps of the $(b+B)$-complexes $\Omh\cinfc(G)$ and $\Omh\cinfc(G_V)$ computing the periodic cyclic homologies $HP_\bullet(\cinfc(G))$ and $HP_\bullet(\cinfc(G_V))$ respectively. We follow the proof of Proposition \ref{plocalg} and introduce infinite row and column matrices $u=(u_i)_{i\in I}$ and $v=(v_i)_{i\in I}$ given by $u_i=c_i\beta_i^{-1}$ and $v_i=\beta_i c_i$. Then $uv=1$, and $vu$ is an idempotent infinite matrix. Moreover $\rho(a)=v a u$ for all $a\in\cinfc(G)$ by definition. The invertible matrices
$$
W = \left(\begin{matrix} 0 & -u \\ v & 1-vu \end{matrix} \right)\ ,\qquad W^{-1} = \left(\begin{matrix} 0 & u \\ -v & 1-vu \end{matrix} \right)
$$
are multipliers of the algebras $M_{\infty}(\cinfc(G))$ and $M_{\infty}(\cinfc(G_V))$. Furthermore, the identity $\rho^1(a)=W^{-1}\rho^0(a)W$ holds for all $a\in\cinfc(G)$, where $\rho^0,\rho^1:\cinfc(G)\to M_{\infty}(\cinfc(G))$ are the homomorphisms
$$
\rho^0(a)= \left(\begin{matrix} a & 0 \\ 0 & 0 \end{matrix} \right) \ ,\qquad \rho^1(a)= \left( \begin{matrix} 0 & 0 \\ 0 & \rho(a) \end{matrix} \right)\ .
$$
Then $\iota_*\circ\tr\rho_*=\tr\rho_*^1$ is chain homotopic to $\tr\rho_*^0 = \id$ on $\Omh\cinfc(G)$. In the same way $\tr\rho_*\circ\iota_*=\tr\rho_*^1$ is chain homotopic to $\tr\rho_*^0=\id$ on $\Omh\cinfc(G_V)$. The isomorphism $HP_\bullet(\cinfc(G))\cong HP_\bullet(\cinfc(G_V))$ follows, as well as the isomorphism in periodic cyclic cohomology. The proof is the same for topological periodic cyclic (co)homology. \hfill\cqfd\\

\begin{lemma}\label{lbott}
Let $G$ be  Lie groupoid. Let $G'$ be the direct product of the pair groupoid $\rr\times\rr$ with $G$. Then one has isomorphisms
\be
HP_\bullet^\top(\cinfc(G))\cong HP_\bullet^\top(\cinfc(G'))\ ,\qquad HP^\bullet_\top(\cinfc(G))\cong HP^\bullet_\top(\cinfc(G'))\ .
\ee
\end{lemma}
{\it Proof:} The convolution algebra of the pair groupoid $\rr\times\rr$ acts as smoothing operators on the Hilbert space of square-integrable functions on $\rr$ with respect to the Lebesgue measure: the action of $k\in \cinfc(\rr\times\rr)$ on a function $f\in L^2(\rr,dx)$ reads
$$
(k\cdot f)(x) = \int_\rr k(x,y)\, f(y)\, dy\ .
$$
Now choose a Hilbert basis $(|e_i\rangle)_{i\in\nn}$ of $L^2(\rr,dx)$ with the following properties: i) each $|e_i\rangle$ is a smooth function with compact support on $\rr$ and ii) the infinite $\nn\times\nn$ matrix with scalar coefficients $k_{ij}=\langle e_i| k |e_j\rangle$ is of rapid decay for all $a\in \cinfc(\rr\times\rr)$. Such a basis can be obtained by modification of the orthonormal basis of Hermite polynomials, in such a way that the unitary matrix passing from the Hermite polynomials to $(|e_i\rangle)$ is of the form $1$ + a matrix with rapid decay. Thus, $\cinfc(\rr\times\rr)$ is represented as a subalgebra of $\Kc$, the algebra of $\nn\times\nn$ matrices with rapid decay. There is natural locally convex topology on $\Kc$ defined by the family of norms
$$
\| k\|_n = \sum_{(i,j)\in\nn\times\nn} (1+i+j)^n|k_{ij}|
$$ 
for all non-negative integers $n$, which turns $\Kc$ into a Fr\'echet algebra containing $\cinfc(\rr\times\rr)$ as a dense subalgebra. In fact the inclusion $\cinfc(\rr\times\rr) \to \Kc$, sending an operator $k$ to the matrix with coeficients $k_{ij}=\langle e_i| k |e_j\rangle$, is continuous with respect to the LF topology on $\cinfc(\rr\times\rr)$. Remark also that the convolution algebra of the product groupoid $G'=(\rr\times\rr)\times G$ is isomorphic to $\cinfc(\rr\times\rr)\hotimes \cinfc(G)$, where $\hotimes$ is the inductive tensor product of LF spaces. From this observation one has a continuous inclusion $\iota: \cinfc(G) \to \cinfc(G')$ defined by $\iota(a) = |e_0\rangle \langle e_0| \otimes a$ for all $a\in\cinfc(G)$, where $|e_0\rangle \langle e_0|\in \cinfc(\rr\times\rr)$ is the orthogonal projector associated to the vector $|e_0\rangle \in L^2(\rr,dx)$. The extension of $\iota$ to differential forms yields a morphism of (topological) cyclic bicomplexes 
$$
\iota_*:\Om_\top\cinfc(G)\to\Om_\top\cinfc(G')\ .
$$
Now write $\Kc(\cinfc(G'))$ for the completed tensor product $\Kc\hotimes\cinfc(G')$. Then $\Kc(\cinfc(G'))$ is a completion of the algebra of infinite matrices $M_{\infty}(\cinfc(G'))$. We construct a continuous homomorphism $\rho:\cinfc(G')\to \Kc(\cinfc(G'))$ by setting the $(i,j)$ entry of the matrix $\rho(b)$ equal to $\rho(b)_{ij} = |e_0\rangle \langle e_0| \otimes \langle e_i| b |e_j\rangle$ for all $b\in\cinfc(G')$. Remark that $\langle e_i| b |e_j\rangle \in \cinfc(G)$, so that $\rho(b)_{ij}$ lies in the image of $\iota$. Hence the extension of $\rho$ to differential forms composed with the matrix trace $\tr:\Kc\to\cc$ and the operator trace $\tr':\cinfc(\rr\times\rr)\to\cc$ gives rise to a morphism of cyclic bicomplexes
$$
\tr'\tr\rho_*: \Om_\top\cinfc(G') \to \Om_\top\cinfc(G)\ .
$$
It remains to show that $\iota_*\circ\tr'\tr\rho_*$ and $\tr'\tr\rho_*\circ\iota_*$ are chain-homotopic to the identity maps on the $(b+B)$-complexes $\Omh\cinfc(G')$ and $\Omh\cinfc(G)$ respectively. For any $a\in\cinfc(G)$ one has $\rho(\iota(a))_{ij}=0$ if $(i,j)\neq (0,0)$ and $\rho(\iota(a))_{00}=\iota(a)$. Hence the ity $\tr'(|e_0\rangle\langle e_0|)=1$ gives $\tr'\tr\rho_*\circ\iota_*=\id$ on $\Omh\cinfc(G)$. In the converse direction, a simple computation yields the equality $\iota_*\circ\tr'\tr\rho_*=\tr\rho_*$, hence it is sufficient to show that $\tr\rho_*$ is chain-homotopic to the identity map on $\Omh\cinfc(G')$. Let $u=(u_i)_{i\in\nn}$ be the infinite row matrix with entries $u_i=|e_i\rangle\langle e_0| \in \cinfc(\rr\times\rr)$, and $v=(v_i)_{i\in\nn}$ be the infinite column matrix with entries $v_i=|e_0\rangle\langle e_i| \in \cinfc(\rr\times\rr)$. We regard $\cinfc(\rr\times\rr)$ as an algebra of multipliers of $\cinfc(G')$. Then $uv=1$ which allows to define as usual the invertible matrices
$$
W = \left(\begin{matrix} 0 & -u \\ v & 1-vu \end{matrix} \right)\ ,\qquad W^{-1} = \left(\begin{matrix} 0 & u \\ -v & 1-vu \end{matrix} \right)\ .
$$
By definition $\rho(b)=vbu$ for all $b\in\cinfc(G')$, hence the homomorphisms $\rho^0,\rho^1:\cinfc(G')\to \Kc(\cinfc(G'))$ given for all $b\in\cinfc(G')$ by 
$$
\rho^0(b)= \left(\begin{matrix} b & 0 \\ 0 & 0 \end{matrix} \right) \ ,\qquad \rho^1(b)= \left( \begin{matrix} 0 & 0 \\ 0 & \rho(b) \end{matrix} \right)\ ,
$$
are conjugate under the adjoint action of $W$. The elements $\rho^0(b)$, $W^{-1}\rho^0(b)$, $\rho^0(b)W$ and $W^{-1}\rho^0(b)W=\rho^1(b)$ are all in $\Kc(\cinfc(G'))$, hence the classical argument using rotation matrices shows that $\rho^0$ and $\rho^1$ are stably homotopic (notice however that $W$ and $W^{-1}$ are not multipliers of the whole algebra $\Kc(\cinfc(G'))$). Consequently $\tr\rho_*$ is chain-homotopic to the identity map on $\Omh\cinfc(G')$ as wanted. The isomorphisms $HP_\bullet^\top(\cinfc(G)) \cong HP_\bullet^\top(\cinfc(G'))$ and $HP^\bullet_\top(\cinfc(G))\cong HP^\bullet_\top(\cinfc(G'))$ follow. \hfill\cqfd\\

\begin{proposition}[Morita invariance]\label{pmor}
Let $G_1\rightrightarrows B_1$ and $G_2\rightrightarrows B_2$ be Lie groupoids and let $B_1\stackrel{\pi_1}{\longleftarrow} M \stackrel{\pi_2}{\longrightarrow} B_2$ be a Morita equivalence. Then one has isomorphisms in topological periodic cyclic (co)homology
\be 
HP^\top_\bullet (\cinfc(G_1)) \cong HP^\top_\bullet (\cinfc(G_2)) \ ,\quad HP_\top^\bullet (\cinfc(G_1)) \cong HP_\top^\bullet (\cinfc(G_2)) \ .
\ee
Moreover, if the submersions $\pi_1$ and $\pi_2$ are \emph{\'etale}, the isomorphisms also hold in algebraic periodic cyclic (co)homology.
\end{proposition}
{\it Proof:} It sufficies to show that, given a groupoid $G\rightrightarrows B$ and a surjective submersion $\pi:M\to B$, the periodic cyclic (co)homologies of $\cinfc(G)$ and $\cinfc(\pi^*G)$ are isomorphic. Let $n$ be the dimension of the fibers of $\pi$, and denote by $\pi_0:B\times\rr^n\to B$ the projection onto the first factor. The manifolds $B\times\rr^n$ and $M$ have the same dimension. Let $U=M\coprod (B\times\rr^n)$ be their disjoint union. $\pi$ and $\pi_0$ thus give a surjective submersion $\si:U\to B$. Then $M$ and $B\times\rr^n$ are open subsets of $U$ intersecting each orbit of the groupoid $\si^*G$. Moreover $\pi^*G$ is the restriction groupoid of $\si^*G$ to the subset $M$, and $\pi_0^*G$ is the restriction groupoid of $\si^*G$ to the subset $B\times\rr^n$. By Lemma \ref{lgeorges}, the periodic cyclic (co)homologies of $\cinfc(\pi^*G)$, $\cinfc(\pi_0^*G)$ and $\cinfc(\si^*G)$ are isomorphic, in the algebraic as well as topological setting. If $n=0$, that is when $\pi$ is \'etale, $\pi_0^*G=G$ implies that $\cinfc(G)$ and $\cinfc(\pi^*G)$ have the same \emph{algebraic} periodic cyclic (co)homology. If $n$ is arbitrary, then $\pi_0^*G$ is the direct product of $G$ with the pair groupoid $\rr^n\times\rr^n$. By virtue of Lemma \ref{lbott}, the \emph{topological} periodic cyclic (co)homologies of $\cinfc(G)$, $\cinfc(\pi_0^*G)$ and $\cinfc(\pi^*G)$ coincide. \hfill\cqfd\\

\section{Excision}\label{sexci}

A convenient way to calculate excision in periodic cyclic cohomology is provided by the Cuntz-Quillen formalism \cite{CQ95,CQ97}. We recall that the cyclic homology of an associative algebra $\Ac$ can be entirely recovered from cyclic $0$-cocycles (traces) and cyclic $1$-cocycles over suitable extensions $0\to\Jc\to\Rc\to\Ac\to 0$ of $\Ac$. The basic ingredient is the $X$-complex of $\Rc$,
\be
X(\Rc)\ :\ \Rc\rightleftarrows \Om^1\Rc_\nat\ ,
\ee
where $\Om^1\Rc_\nat = \Om^1\Rc/[\Rc,\Om^1\Rc]$. For all elements $x,y\in \Rc$ we write $\nat x\dd y$ for the class of the one-form $x\dd y \mod [\Rc,\Om^1\Rc]$. The boundary map $\d_0:\Rc \to \Om^1\Rc_\nat$ is the non-commutative differential $x\mapsto \nat\dd x$, while the boundary map $\d_1:\Om^1\Rc_\nat \to \Rc$ is the commutator $\nat x\dd y \mapsto [x,y]$ induced by the Hochschild boundary on $\Om^1\Rc$. One has $\d_1\circ\d_0=0=\d_0\circ\d_1$ hence $X(\Rc)$ is a $\zz_2$-graded complex, with even part $\Rc$ and odd part $\Om^1 \Rc_\nat$. If $\Jc\subset \Rc$ is a two-sided ideal, Cuntz and Quillen define a decreasing filtration of $X(\Rc)$ by the subcomplexes $F^n_{\Jc}X(\Rc)$, $n\in\zz$, as follows:
\beq
F^{2n}_{\Jc}X(\Rc)  &:& \Jc^{n+1} + [\Jc^n,\Rc] \rightleftarrows \nat\Jc^n_{(+)}\dd\Rc \\
F^{2n+1}_{\Jc}X(\Rc) &:& \Jc^{n+1} \rightleftarrows \nat \big(\Jc^{n+1}_{(+)}\dd\Rc + \Jc^n_{(+)}\dd\Jc \big)
\eeq
where $\Jc^n=\Rc$ and $\Jc^n_{(+)}=\Rc^+$ for $n\leq 0$. In particular $F^n_{\Jc}X(\Rc) = X(\Rc)$ whenever $n<0$. The $\Jc$-adic completions of $\Rc$ and $X(\Rc)$ are respectively a pro-algebra and a pro-complex given by the projective limits
\be
\Rch = \varprojlim_n \Rc/\Jc^n\ ,\qquad X(\Rch) = \varprojlim_n X(\Rc)/F^n_{\Jc}X(\Rc)\ .
\ee
A cocycle of pro-complex $\tau : X(\Rch) \to \cc$ is exactly a cocycle over $X(\Rc)$ vanishing on the sub complex $F^n_{\Jc}X(\Rc)$ for some $n$. Thus, a cocycle of even degree is a trace on $\Rc$ vanishing on the large powers of the ideal $\Jc$. Similarly, a cocycle of odd degree is a cyclic $1$-cocycle over $\Rc$ vanishing whenever one of its arguments lies in $\Jc^n$ for some $n$. The link with the cyclic homology of the quotient algebra $\Ac = \Rc/\Jc$ shows up in the case of the universal free extension $\Rc=T\Ac$ corresponding to the non-unital tensor algebra of $\Ac$:
\be
T\Ac = (\Ac) \oplus (\Ac\otimes\Ac) \oplus (\Ac\otimes\Ac\otimes\Ac) \oplus \ldots
\ee
The product on $T\Ac$ is the tensor product, and by definition the two-sided ideal $\Jc=J\Ac$ is the kernel of the multiplication homomorphism $m:T\Ac\to\Ac$, $a_1\otimes\ldots\otimes a_n \mapsto a_1\ldots a_n$. One checks that $J\Ac$ is generated by the inhomogeneous elements of the form $a_1a_2 - a_1\otimes a_2$. As a $\zz_2$-graded vector space, $X(T\Ac)$ is isomorphic to the space of non-commutative differential forms $\Om\Ac$.  More precisely $T\Ac$ is isomorphic to the space $\Om^+\Ac$ of differential forms of even degree over $\Ac$, whereas $\Om^1T\Ac_\nat\cong T\Ac^+\otimes\Ac$ is isomorphic to the space $\Om^-\Ac$ of differential forms of odd degree. These isomorphisms are given by
\beq
&&a_0da_1\ldots da_{2n}\in \Om^{2n}\Ac \leftrightarrow a_0\otimes \om(a_1,a_2)\otimes\ldots \otimes \om(a_{2n-1},a_{2n}) \in T\Ac\ ,\non\\
&&a_0da_1\ldots da_{2n+1} \in \Om^{2n+1}\Ac \leftrightarrow \nat(a_0\otimes \om(a_1,a_2)\otimes\ldots \dd a_{2n+1}) \in \Om^1T\Ac_\nat\ ,\non
\eeq
where $\om(a_1,a_2)=a_1a_2-a_1\otimes a_2$. The tensor product on $T\Ac$ corresponds to the Fedosov product on $\Om^+\Ac$, given by $\om_1\odot\om_2=\om_1\om_2-d\om_1 d\om_2$ for all differential forms $\om_1,\om_2$ of even degree. The ideal $J\Ac$ corresponds to the space of all forms of even degree $\geq 2$. Cuntz and Quillen show in \cite{CQ95} that the linear isomorphism $\Om\Ac\cong X(T\Ac)$ is, up to a rescaling factor, a quasi-isomorphism between the cyclic bicomplex of $\Ac$ and the $X$-complex of $T\Ac$ endowed with the filtration $F^n_{J\Ac}X(T\Ac)$. In particular the periodic cyclic homology of $\Ac$ is the homology of the pro-complex $X(\Th\Ac)$. More generaly $X(\Rch)$ computes the periodic cyclic homology of $\Ac$ provided that $0\to\Jc\to\Rc\to\Ac\to 0$ is a \emph{quasi-free} extension, see \cite{CQ95}.\\

Let $(E):0\to\Bc\to\Ec\to\Ac\to 0$ be an arbitrary extension. As proved by Cuntz and Quillen in \cite{CQ97}, this extension gives rise to an excision map in periodic cyclic cohomology
\be 
E^*\ :\ HP^\bullet(\Bc)\to HP^{\bullet+1}(\Ac)\ .
\ee 
We shall explain how to compute it, following \cite{P8}. The universal property of the tensor algebra allows to lift the homomorphism $\Ec\to\Ac$ to an homomorphism $T\Ec\to T\Ac$ sending the ideal $J\Ec$ to $J\Ac$. We obtain in this way a commutative diagram where all rows and columns are extensions:
\be
\vcenter{\xymatrix{
 & 0 \ar[d] & 0 \ar[d] & 0 \ar[d] &  \\
0 \ar [r] & \Jc \ar[r] \ar[d] & J\Ec \ar[r] \ar[d] & J\Ac  \ar[r] \ar[d] & 0  \\
0 \ar [r] & \Rc \ar[r] \ar[d] & T\Ec \ar[r] \ar[d] & T\Ac  \ar[r] \ar[d] & 0  \\
0 \ar [r] & \Bc \ar[r] \ar[d] & \Ec \ar[r] \ar[d] & \Ac  \ar[r] \ar[d] & 0  \\
 & 0 & 0 & 0 &  }}\label{diag}
\ee
By construction $\Rc$ is the kernel of the homomorphism $T\Ec\to T\Ac$, and $\Jc$ is the kernel of $J\Ec\to J\Bc$. The above diagram allows to calculate the excision map $HP^\bullet(\Bc)\to HP^{\bullet+1}(\Ac)$, once the relevant cyclic cohomology classes of $\Bc$ are represented in a suitable form. To this end we note that (\ref{diag}) allows to define two relevant filtrations of the complex $X(T\Ec)$. The first filtration is induced by the ideal $J\Ec$. We denote the corresponding $J\Ec$-adic completions with a hat:
\be 
\Th\Ec = \varprojlim_k T\Ec/(J\Ec)^k\ ,\qquad X(\Th\Ec) = \varprojlim_k X(T\Ec)/F^k_{J\Ec}X(T\Ec)\ .
\ee
The complex $X(\Th\Ec)$ computes the periodic cyclic homology of $\Ec$. The second filtration is induced by the ideal $J\Ec+\Rc$. We denote the corresponding completions with a tilde. Then we can write
\be 
\widetilde{T}\Ec = \varprojlim_n T\Ec/(J\Ec+\Rc)^n\ ,\qquad X(\widetilde{T}\Ec) = \varprojlim_n X(\Th\Ec)/F^n_{\Rc}(\Th\Ec)
\ee
where $F^n_{\Rc}X(\Th\Ec) = \varprojlim_k F^n_{J\Ec}X(T\Ec)/(F^n_{J\Ec}X(T\Ec) \cap F^k_{\Rc}X(T\Ec))$ is a subcomplex of $X(\Th\Ec)$.

\begin{lemma}[\cite{P8}]
The $\zz_2$-graded complex of cochains over $F^n_{\Rc}X(\Th\Ec)$ computes the periodic cyclic cohomology $HP^\bullet(\Bc)$ for all $n\geq 1$. \hfill\cqfd
\end{lemma}
Hence any periodic cyclic cohomology class over $\Bc$ can be represented by a cocycle $\tau: F^n_{\Rc}X(\Th\Ec)\to\cc$ for some $n\geq 1$. We will see below interesting examples of such cocycles in the case of Lie groupoids. Let us now describe in full generality the excision map associated to the extension $0\to\Bc\to\Ec\to\Ac\to 0$, sending a cocycle $\tau$ over $F^n_{\Rc}X(\Th\Ec)$ to a cyclic cocycle over $\Ac$. This is explained in \cite{P8} \S 2. We first choose any extension of $\tau$ to a \emph{linear} map
\be 
\tau_R: X(\Th\Ec)\to \cc\ .
\ee
We call $\tau_R$ a renormalization of $\tau$. Of course $\tau_R$ is generally not a cocycle in $\hom(X(\Th\Ec),\cc)$. However its composite $\tau_R\d$ with the $X$-complex boundary map $\d$ is a cocycle vanishing on the subcomplex $F^n_{\Rc}X(\Th\Ec)$ by construction. Hence $\tau_R\d$ defines a cocycle in $\hom(X(\widetilde{T}\Ec),\cc)$. Then observe that the diagonal of (\ref{diag}) leads to an extension of $\Ac$ by the algebra $T\Ec$, with kernel the ideal $J\Ec+\Rc$. Choose any linear splitting $\si:\Ac\to T\Ec$ of this extension. The universal property of the tensor algebra $T\Ac$ allows to lift the linear map $\si$ to a homomorphism of algebras $\si_*:T\Ac\to T\Ec$ respecting the ideals:
\be
\vcenter{\xymatrix{0 \ar[r] & J\Ac \ar[r] \ar[d]_{\si_*} & T\Ac \ar[r] \ar[d]_{\si_*} & \Ac \ar[r] \ar@{.>}[dl]^{\si} \ar@{=}[d] & 0 \\
0\ar[r] & J\Ec + \Rc \ar[r]  & T\Ec \ar[r] & \Ac \ar[r]   & 0 }}
\ee
Explicitly $\si_*(a_1\otimes\ldots\otimes a_n) = \si(a_1)\ldots\si(a_n)$ in $T\Ec$. Since $\si_*$ respects the ideals, it extends to a homomorphism $\Th\Ac\to\widetilde{T}\Ec$. This in turn induces a chain map still denoted $\si_*:X(\Th\Ac)\to X(\widetilde{T}\Ec)$.

\begin{proposition}[\cite{P8}]
The excision map $HP^\bullet(\Bc)\to HP^{\bullet+1}(\Ac)$ associated to the extension $0\to\Bc\to\Ec\to\Ac\to 0$ is realized by sending a cocycle $\tau \in \hom(F^n_{T(\Bc:\Ec)}X(\Th\Ec),\cc)$ to the cocycle $\tau_R\d\circ\si_* \in\hom(X(\Th\Ac),\cc)$ for any choice of renormalization $\tau_R$ and linear splitting $\si:\Ac\to T\Ec$. \hfill\cqfd
\end{proposition}

We now apply this formalism to the pseudodifferential extensions obtained for groupoid actions in section \ref{sgrou}. Let  $G\rightrightarrows B$ be a Lie groupoid and denote by 
$$
\Cc = \cinfc(B)\cp G
$$ 
the corresponding convolution algebra. Then any cyclic cohomology class $[\varphi] \in HP^\bullet(\Cc)$ is represented by a cocycle $\varphi\in \hom(X(\Th\Cc),\cc)$ where $\Th\Cc$ is the $J\Cc$-adic completion of the tensor algebra $T\Cc$. Taking the locally convex topology of $\Cc$ into account, we find that $HP^\bullet_\top(\Cc)$ is the cohomology of the complex $\hom(X(\Th\Cc_\top),\cc)$ of \emph{continuous} and \emph{bounded} cochains. A continuous bounded cochain of even degree is a linear map $\varphi:T\Cc\to \cc$ given by a family of distributions $\varphi^+_n\in C^{-\infty}(G^n)$, $n\geq 1$, with singularity order bounded uniformly in $n$, such that
\be 
\varphi(f_1\otimes\ldots\otimes f_n) = \int_{G^n} \varphi^+_n(g_1,\ldots,g_n)\, f_1(g_1)\ldots f_n(g_n)
\ee
for all $f_i \in \Cc$, and $\varphi$ has to vanish on the large powers of $J\Cc$. In the same way, a continuous cochain of odd degree is a linear map $\varphi: \Om^1T\Cc_\nat \cong T\Cc^+\otimes \Cc \to \cc$ given by a family of distributions $\varphi^-_n\in C^{-\infty}(G^n)$, $n\geq 1$, with singularity order bounded uniformly in $n$, such that
\be 
\varphi(\nat(f_1\otimes\ldots\otimes f_{n-1}\dd f_n)) = \int_{G^n} \varphi^-_n(g_1,\ldots,g_{n-1}|g_n) f_1(g_1)\ldots f_n(g_n)
\ee
for all $f_i \in \Cc$, and $\varphi$ has to vanish whenever its argument lies in a large power of $J\Cc$. In view of the quasi-isomorphism of complexes $\Omh\Cc \cong X(\Th\Cc)$, the localization proposition \ref{ploctop} shows that any topological cyclic cohomology class can be represented by an $X$-complex cocycle $\varphi=(\varphi^\pm_n)_{n\geq 1}$ such that each distribution $\varphi^\pm_n$ has a support contained in an arbitrary small neighborhood of the set of loops $I^{(n)}\subset G^n$. \\

Now let $\pi:M\to B$ be a $G$-equivariant submersion. We defined in section \ref{sgrou} the crossed product algebra of $G$ with the sections of the bundle $\CL^0_c(M)\to B$ of vertical compactly supported pseudodifferential operators of order $\leq 0$. This leads to an extension $0\to\Bc\to\Ec\to\Ac\to 0$ of the convolution algebra of the groupoid $S^*M_B\rtimes G$, where
$$
\Bc = \cinfc(B,\CL^{-1}_c(M))\cp G\ ,\ \Ec = \cinfc(B,\CL^0_c(M))\cp G\ ,\ \Ac = \cinfc(S^*_\pi M)\rtimes G \ .
$$
Recall that the ideal $\cinfc(B,\L^{-\infty}_c(M))\cp G \subset \Bc$ of smoothing operators is canonically isomorphic to the smooth convolution algebra of the pullback groupoid $\pi^*G\rightrightarrows M$. Hence, any topological cyclic cohomology class $[\varphi]\in HP^\bullet_\top(\Cc)$ corresponds to a class in $HP^\bullet_\top(\cinfc(B,\L^{-\infty}_c(M))\cp G)$ by virtue of the Morita equivalence of groupoids $G\sim \pi^*G$. Our first goal is to show that $[\varphi]$ can be lifted to a class $[\tau_\varphi] \in HP^\bullet(\Bc)$. This requires the following notion of connection.

\begin{definition}\label{dcon}
A \emph{generalized connection} on a submersion $\pi:M\to B$ is a function $h$ defined on a neighborhood of the diagonal in $B\times B$, sending any pair of neighboring points $(b_1,b_2)$ to a linear operator $h(b_1,b_2): \cinfc(M_{b_2}) \to \cinfc(M_{b_1})$ which may be decomposed as a locally finite sum
\be
h(b_1,b_2)= \sum_{i} h_i(b_1,b_2) \circ  f_i(b_1,b_2) 
\ee
where:\\
$\bullet$ $h_i({b_1},{b_2})$ is a local diffeomorphism from an open subset of $M_{b_2}$ to an open subset of $M_{b_1}$, depending smoothly on $(b_1,b_2)$, such that $h_i(b,b)=\Id$ for all $b\in B$;\\
$\bullet$ $f_i(b_1,b_2) $ is a smooth function on $M_{b_2}$ with support contained in the domain of $h_i(b_1,b_2)$, depending smoothly on $(b_1,b_2)$, such that $\sum_i f_i(b,b)=1$ for all $b\in B$. The action of $f_i(b_1,b_2)$ on $\cinfc(M_{b_2})$ is by pointwise multiplication.
\end{definition}
A submersion always has a generalized connection. Indeed let $n=\dim(M/B)$ be the dimension of the fibers of $\pi:M\to B$ and consider the trivial submersion $\pi_0:B\times\rr^n \to B$ as in the proof of Proposition \ref{pmor}. Then we can find a locally finite open covering $(U_i)_{i\in I}$ of $M$ together with local diffeomorphisms $\beta_i:U_i\to V_i\subset B\times\rr^n$ compatible with the projections, i.e. $\pi_0\circ\beta_i = \pi\vert_{U_i}$ for all $i\in I$. Let $(c_i)_{i\in I}$, $\sum_i c_i^2=1$, be a smooth partition of unity relative to the covering $(U_i)$. For all $b\in B$ let $\beta_i(b) : U_i\cap M_b\to V_i\cap \{b\}\times \rr^n$ and $c_i(b)\in\cinf(M_b)$ be respectively the restriction of the diffeomorphism $\beta_i$ and the function $c_i$ to the fiber $M_b$. Since the fibers of $B\times \rr^n$ are canonically diffeomorphic, we can view $\beta_i(b_1)^{-1}\beta_i(b_2)$ as a local diffeomorphism from $M_{b_2}$ to $M_{b_1}$ for all pairs $(b_1,b_2)\in B\times B$. Then 
\be
h(b_1,b_2) = \sum_{i\in I} c_i(b_1)\circ \beta_i(b_1)^{-1}\beta_i(b_2) \circ c_i(b_2)\ ,
\ee
where $c_i(b_1)$ and $c_i(b_2)$ act by pointwise multiplication on the vector spaces $\cinfc(M_{b_1})$ and $\cinfc(M_{b_2})$ respectively, defines a generalized connection on $M$, with local diffeomorphisms $h_i(b_1,b_2) = \beta_i(b_1)^{-1}\beta_i(b_2)$ and smooth functions $f_i(b_1,b_2) = (c_i(b_1)\circ h_i(b_1,b_2)) \, c_i(b_2)$. In this example the function $h$ is defined on the entire product $B\times B$ and not only on a neighborhood of its diagonal.

\begin{remark}\label{rcon}
\textup{In fact one can always find a generalized connection where the local diffeomorphisms $h_i$ agree pairwise on their common domain in a neighborhood of the diagonal. A possible construction goes as follows. Choose a Riemannian metric on $B$ and a horizontal distribution on $M$, that is, a subbundle $H'$ of the tangent bundle $TM$ complementary the vertical tangent bundle $\ker(\pi_*:TM\to TB)$. That is, one has a decomposition $TM = H'\oplus \ker(\pi_*)$. Then on a suitable open subset the diffeomorphism $h_i({b_1},{b_2})$ is obtained by lifting the geodesic between ${b_2}$ and ${b_1}$ according to the horizontal paths determined by $H'$.}
\end{remark}
We now come back to the situation of the groupoid $G\rightrightarrows B$ acting on the submersion $\pi:M\to B$, and fix any choice of generalized connection $h$. Let $e_1,\ldots, e_n$ be elements of the crossed product $\Ec= \cinfc(B,\CL^0_c(M))\cp G$, such that each $e_i$ belongs to the subspace $\cinfc(B,\CL^{-m_i}_c(M))\cp G$ with $m_1+\ldots +m_n >\dim(M/B)+k$ for a given integer $k$. We can map the tensor $e_1 \otimes \ldots \otimes e_n \in T\Ec$ to a compactly supported function $\Tr^{h}_{e_1 , \ldots , e_n}$ of class $C^k$ on an appropriate neighborhood $V$ of the set of loops $I^{(n)}\subset G^n$. This function is  defined by evaluation on any point $(g_1,\ldots,g_n)\in V$ by
\beq
\lefteqn{\Tr^{h}_{e_1, \ldots , e_n}(g_1,\ldots, g_n) = }  \label{toto} \\ 
&& \qquad \Tr\big(e_1(g_1)\, U_{g_1}\, {h(s(g_1),r(g_2))}  \ldots e_n(g_n)\, U_{g_n}\, {h(s(g_n),r(g_1))}\big) \non
\eeq
where $\Tr$ is the ordinary trace of operators acting on the space of scalar functions on the manifold $M_{r(g_1)}$. In order to show that this expression makes sense, recall that $e_i(g_i) \in \CL^{-m_i}_c(M_{r(g_i)})$ and $U_{g_i}$ is the diffeomorphism from $M_{s(g_i)}$ to  $M_{r(g_i)}$ defined by the action of $g_i$ on $M$. Since $h(s(g_i),r(g_{i+1})) $ is a sum of diffeomorphisms composed with pointwise multiplication by smooth functions, the product $e_i(g_i)\, U_{g_i}\, h(s(g_i),r(g_{i+1})) $ is a compactly supported operator carrying smooth functions on $M_{r(g_{i+1})}$ to smooth functions on $M_{r(g_i)}$. Hence the product under the trace in (\ref{toto}) is a sum of pseudodifferential operators in $\CL^{-m}_c(M_{r(g_1)})$ composed with diffeomorphisms of $M_{r(g_1)}$, with $m=m_1+\ldots + m_n$. Its partial derivatives of order $k$ with respect to the variables $g_i$ yield a sum of operators in $\CL^{-m+k}_c(M_{r(g_1)})$ composed with diffeomorphisms, which remains in the domain of the trace. Therefore $\Tr^{h}_{e_1, \ldots , e_n}$ is a function of class $C^k$. 

\begin{lemma}\label{ltau}
As above let $\Cc = \cinfc(B)\cp G$, $\Ec = \cinfc(B,\CL^0_c(M))\cp G$, $\Ac = \cinfc(S^*_\pi M)\cp G$ and $\Rc=\ker(T\Ec\to T\Ac)$. Choose any generalized connection $h$ on the submersion $\pi:M\to B$. Then the map sending a bounded continuous cochain $\varphi\in \hom(X(\Th\Cc_\top),\cc)$ to the cochain $\tau_\varphi\in \hom(F^m_{\Rc}X(\Th\Ec),\cc)$ defined by
\beq
&&\tau_\varphi(e_1\otimes\ldots\otimes e_n) = \int_{G^n} \varphi^+_n(g_1,\ldots,g_n)\, \Tr^h_{e_1,\ldots,e_n}(g_1,\ldots,g_n)\ , \non\\
&&\tau_\varphi(\nat(e_1\otimes\ldots\otimes e_{n-1}\dd e_n)) = \int_{G^n} \varphi^-_n(g_1,\ldots,g_{n-1}|g_n)\, \Tr^h_{e_1,\ldots,e_n}(g_1,\ldots,g_n) \non
\eeq
for all $e_i\in \Ec$, $n\geq 1$, is a morphism of complexes provided that $m$ is sufficiently large. The induced map
\be 
\tau_* \ :\ HP_\top^{\bullet}(\cinfc(B)\cp G) \to HP^\bullet(\cinfc(B,\CL_c^{-1}(M))\cp G)
\ee
does not depend on the choice of connection $h$.
\end{lemma}
{\it Proof:} The distributions $\varphi^\pm$ being of finite singularity order, $\tau_\varphi$ is well-defined provided that the function $\Tr^{h}_{e_1, \ldots , e_n}$ is regular enough, that is, $m$ sufficiently large. Since $h=\Id$ on the diagonal of $B\times B$, one easily checks that the trace map $e_1\otimes\ldots\otimes e_n \mapsto \Tr^{h}_{e_1, \ldots , e_n}$ commutes with all operators on the cyclic bicomplexes of $\Ec$ and $\Cc$. Therefore $\tau_\varphi$ is a cocycle. The fact that two different choices of generalized connections give cohomologous cocycles is a consequence of a classical transgression formula. \cqfd\\

\begin{proposition}\label{ptoto}
For any Lie groupoid $G\rightrightarrows B$ and any $G$-equivariant surjective submersion $\pi:M\to B$, one has a commutative diagram
\be 
\vcenter{\xymatrix{
HP_\top^{\bullet}(\cinfc(M)\cp \pi^*G) \ar[r] & HP^\bullet(\cinfc(M)\cp \pi^*G) \\
HP_\top^{\bullet}(\cinfc(B)\cp G) \ar@{=}[u] \ar[r]^{\tau_*\quad} & HP^\bullet(\cinfc(B,\CL_c^{-1}(M))\cp G) \ar[u] }}
\ee
where the left equality is the canonical isomorphism accounting for the Morita equivalence of groupoids $G\sim \pi^*G$, and the right vertical arrow is the restriction morphism induced by the inclusion of the ideal $\cinfc(M)\cp\pi^*G\cong \cinfc(B,\L^{-\infty}_c(M))\cp G$ into $\cinfc(B,\CL^{-1}_c(M))\cp G$.
\end{proposition}
{\it Proof:} If $e_1, \ldots , e_n$ are smoothing operators, the trace map $e_1\otimes\ldots\otimes e_n \mapsto \Tr^{h}_{e_1, \ldots , e_n}$ is precisely the realization of the Morita equivalence $G\sim \pi^*G$ at the level of topological cyclic cohomology. \hfill\cqfd\\

As remarked in section \ref{sgrou}, there is a canonical morphism between the two extensions $(E_0)$ and $(E)$,
$$
\vcenter{\xymatrix{
0 \ar [r] & \cinfc(M)\cp \pi^*G \ar[r] \ar[d] & \cinfc(B,\CL^0_c(M))\cp G \ar[r] \ar@{=}[d] & \cinfc(B,\CS^0_{c}(M))\cp G  \ar[r] \ar[d] & 0  \\
0 \ar [r] & \cinfc(B,\CL^{-1}_{c}(M))\cp G \ar[r] & \cinfc(B,\CL^0_{c}(M))\cp G \ar[r] & \cinfc(S^*_\pi M)\cp G  \ar[r] & 0 }}
$$
where the left vertical arrow is the inclusion homomorphism and the right vertical arrow is the leading symbol homomorphism. By naturality, the respective excision maps $E^*_0$ and $E^*$ are compatible through the induced morphisms in periodic cyclic cohomology, and this combined with Proposition \ref{ptoto} leads to a commutative diagram
$$
\vcenter{\xymatrix{
HP^\bullet(\cinfc(M)\cp\pi^*G) \ar[r]^{\!\!\! E^*_0\quad} & HP^{\bullet+1}(\cinfc(B,\CS^0_{c}(M))\cp G) \\
HP^\bullet(\cinfc(B,\CL_{c}^{-1}(M))\cp G) \ar[u] \ar[r]^{\quad E^*} & HP^{\bullet+1}(\cinfc(S^*_\pi M)\rtimes G) \ar[u] \\
HP_\top^{\bullet}(\cinfc(B)\cp G) \ar[u]^{\tau_*} & }}
$$
As a consequence of Morita invariance, the excision map $E^*_0$ restricted to the image of topological cyclic cohomology of $\cinfc(M)\cp \pi^*G$ thus factors through the cyclic cohomology of the leading symbol algebra $\cinfc(S^*_\pi M\cp G)$. Hence all the relevant information is carried by the excision map $E^*$. 

\begin{theorem}\label{ttop}
Let $G\rightrightarrows B$ be a Lie groupoid and $\pi:M\to B$ a $G$-equivariant surjective submersion. Then one has a commutative diagram
\be 
\vcenter{\xymatrix{
HP^\bullet(\cinfc(B,\CL_{c}^{-1}(M))\cp G) \ar[r]^{\quad E^*} & HP^{\bullet+1}(\cinfc(S^*_\pi M)\rtimes G)  \\
HP_\top^{\bullet}(\cinfc(B)\cp G) \ar[u]^{\tau_*} \ar[r]^{\pi^!_G\quad} & HP_\top^{\bullet+1}(\cinfc(S^*_\pi M)\rtimes G) \ar[u] }}
\ee
\end{theorem}
{\it Proof:} It suffices to remark that for any cocycle $\varphi\in \hom(X(\Th\Cc_\top),\cc)$, the image $\tau_\varphi\in \hom(F^m_{\Rc}X(\Th\Ec),\cc)$ can be renormalized in a continuous and bounded way. For example, by inserting in the trace map (\ref{toto}) a projection operator onto pseudodifferential operators of sufficiently low order.   \cqfd\\

For completeness let us recall the link between the excision map $E^*$ in periodic cyclic cohomology and the $K$-theoretic index map \cite{Ni,P8}. The Chern-Connes pairing between an even cyclic cohomology class $[\tau]\in HP^0(\Bc)$, represented by a trace $\tau:\Rc^m\to \cc$, and a $K$-theory class $[e]-[e_0]\in K_0(\Bc)$ is described as follows (see \cite{P8} \S 4). Suppose for simplicity that $[e]$ is represented by a $2\times 2$ matrix idempotent $e\in M_2(\Bc^+)$ satisfying the property $e-e_0\in M_2(\Bc)$, where $e_0=\bigl(\begin{smallmatrix} 1 & 0 \\ 0 & 0 \end{smallmatrix} \bigr)$. Looking at (\ref{diag}) we see that the intersection $\Rc\cap (J\Ec)^l$ is a two-sided ideal in $\Rc$ for any $l\geq 1$, and the quotient algebra $\Rc/(\Rc\cap (J\Ec)^l)$ is a nilpotent extension of $\Bc$. Define the projective limit
\be
\Rch = \varprojlim_l \Rc/(\Rc\cap (J\Ec)^l)\ .
\ee
It is a classical result (\cite{CQ95}) that $e$ can be lifted to an idempotent $\eh\in M_2(\Rch^+)$, such that $\eh-e_0\in M_2(\Rch)$. Two different liftings give rise to the same $K$-theory class of $\Rch$. Since the trace $\tau$ vanishes on $\Rc^m\cap(J\Ec)^l$ for large $l$, it extends to a trace on $\Rch^m$. Therefore the  pairing ($\tr$ denotes the trace of $2\times 2$ matrices over $\cc$)
\be
\langle [\tau] , [e] \rangle = \tau\#\tr((\eh-e_0)^{2n+1})
\ee
is well-defined for $2n+1\geq m$, and one shows that it only depends on the $K$-theory class $[e]$ and the periodic cyclic cohomology class $[\tau]$. The formula still works for idempotents $e\in M_{\infty}(\Bc^+)$, where $\tr$ is now the trace of matrices of arbitrary size over $\cc$. \\
On the other hand, the odd cyclic cohomology class $E^*([\tau])\in HP^1(\Ac)$ is represented by the cyclic 1-cocycle $\tau_R\d\circ\si_*$ over $\Th\Ac$, and can be paired with any algebraic $K$-theory class of degree one $[u]\in K_1(\Ac)$ as follows. Choose an invertible matrix $u\in GL_{\infty}(\Ac)$ representing the class $[u]$. Then any lifting $\uh\in GL_{\infty}(\Th\Ac)$ of $u$ is invertible. In particular if one chooses $\uh=u$ via the canonical linear map $\Ac\to T\Ac\to\Th\Ac$, its inverse is given by the series
\be
\uh^{-1} = \sum_{n=0}^{\infty} u^{-1}\otimes (1-u\otimes u^{-1})^{\otimes n}\ ,
\ee
where $u^{-1}$ is the inverse of $u$ in $GL_{\infty}(\Ac)$. Note that $1-u\otimes u^{-1}$ belongs to the ideal $M_{\infty}(J\Ac)$ hence the series $\uh^{-1}$ is convergent in $GL_{\infty}(\Th\Ac)$. The cup-product of the cyclic cocycle $\tau_R\d\circ\si_*$ with the trace of matrices over $\cc$ yields the Chern-Connes pairing\footnotemark
\be
\langle E^*([\tau]), [u] \rangle = (\tau_R\d\circ\si_*)\#\tr(\uh^{-1},\uh) = \tau_R\#\tr([\si_*(\uh^{-1}),\si_*(\uh)])
\ee
which only depends on the class $[u]\in K_1(\Ac)$ and the periodic cyclic cohomology class $E^*([\tau])\in HP^1(\Ac)$. The main result of \cite{Ni,P8} is that the excision map in cyclic cohomology is adjoint to the index map in algebraic $K$-theory. Let us state this general theorem in the specialized case of groupoids:
\footnotetext{In \cite{P8} we used different normalization conventions for the Chern-Connes pairing, involving the numerical factor $\sqrt{2\pi \i}$ in the odd case, which was dictated by the need of compatibility with the bivariant Chern-Connes character. Since we don't want to discuss these matters here we use simpler conventions.}

\begin{corollary}[\cite{P8}]\label{cind}
Let $G\rightrightarrows B$ be a Lie groupoid with convolution algebra $\Cc = \cinfc(B)\cp G$, and let $\pi:M\to B$ be a $G$-equivariant surjective submersion. Let
$$
(E)\ :\ 0\to\Bc\to\Ec\to\Ac\to 0
$$ 
be the associated extension of the algebra $\Ac=\cinfc(S^*_\pi M)\cp G$ of non-commutative symbols by $\Bc=\cinfc(B,\CL_c^{-1}(M))\cp G$. Then for any $[u]\in K_1(\Ac)$ and any $[\varphi]\in HP_\top^0(\Cc)$, the pairing of the index $\Ind_E([u])\in K_0(\Bc)$ with the cyclic cohomology class $[\tau_\varphi]\in HP^0_\top(\Bc)$ is given by the formula
\be
\langle [\tau_\varphi] , \Ind_E([u])\rangle = \langle E^*([\tau_\varphi]), [u] \rangle
\ee
where $E^*([\tau_\varphi])\in HP^1(\Ac)$ is represented by the above cocycle $(\tau_\varphi)_R\d\circ\si_*$, for any choice of renormalization $(\tau_\varphi)_R$. \cqfd
\end{corollary}

\section{The residue formula}\label{szeta}

As before consider a Lie groupoid $G\rightrightarrows B$ and a $G$-equivariant submersion $\pi:M\to B$. Let $O\subset G$ be an isotropic submanifold invariant under the adjoint action of $G$. We will compute the excision map on the localized cyclic cohomology classes $[\varphi]\in HP^\bullet_\top(\cinfc(B)\cp G)_{[O]}$ by means of a residue formula. This closely follows (and actually generalizes) the construction of \cite{P8}. In order to make everything work, we need to impose some constraints on the structure of the fixed points for $O$. Remark that any $h\in O$ verifies $r(h)=s(h)$ by definition, hence the fiber $M_{r(h)}$ carries an action of $h$ by diffeomorphisms.  

\begin{definition}
Let $\pi: M\to B$ be a $G$-equivariant submersion and $O\subset G$ an isotropic submanifold. We say that the action of $O$ on $M$ is \emph{non-degenerate} if the following holds: \\

\noindent i) For any $h\in O$, the set of fixed points $M^h_{r(h)}$ is a union of isolated submanifolds in $M_{r(h)}$, depending smoothly on $h$;\\
ii) At any point $x\in M_{r(h)}^h$ the tangent space $T_xM_{r(h)}$ in the ambient manifold $M_{r(h)}$ splits as a direct sum
\be
T_xM_{r(h)} = T_xM_{r(h)}^h\oplus N_x^h
\ee
of two subspaces globally invariant by the action of the tangent map $h_*$ associated to the diffeomorphism. We denote $h'$ the restriction of $h_*$ to the normal subspace $N_x^h$;\\
iii) The endomorphism $1-h'$ of $N_x^h$ is non-singular, that is $\det(1-h')\neq 0$ at any point $x\in M_{r(h)}^h$.
\end{definition}
The non-degeneracy condition is automatically satisfied, for example, when $h$ acts isometrically with respect to a Riemannian metric on $M_{r(h)}$. In the latter case the subspace $N_x^h$ is the fiber of the normal bundle, in the Riemannian sense, of the fixed submanifold $M_{r(h)}^h$ at $x$. Note that when $h$ is not an isometry, condition {\it iii)} may definitely fail; this happens for example in the situation of conformal mappings considered in \cite{P7}. We will not cover such situations in this article. Let us now focus on a connected component of the submanifold $M_{r(h)}^h$, say of dimension $r$ and codimension $s$. We fix a local coordinate system $x=(x_1,\ldots,x_r)$ of $M_{r(h)}^h$, and complete it with a normal coordinate system $y=(y_1,\ldots,y_s)$ in a neighborhood of the fixed submanifold with the following properties:
\begin{itemize}
\item $y=(0,\ldots,0)$ on the fixed submanifold $M_{r(h)}^h$;
\item The tangent vectors $\d/\d y_i$, $i=1$,...,$s$, belong to the normal subspace $N_x^h$ at any point $x\in M_{r(h)}^h$.
\end{itemize}
Such a local coordinate system $(x;y)$ on $M_{r(h)}$ will be called \emph{adapted to the fixed submanifold}. The stability of the subspaces $N^h_x$ implies the following important fact: if $h^*x$ denotes the pullback of the coordinate functions $x$ by the diffeomorphism $h$, then the difference $x-h^*x$ is of order $2$ with respect to the variable $y$ near the fixed submanifold $y=0$, while $y-h^*y$ is only of order $1$. This will be used for establishing the properties of zeta-functions.\\

Let $b\in B$ be a point and $Q\in \CL^1(M_b)$ be a properly supported, elliptic, positive and invertible pseudodifferential operator of order one. For example, we can take $Q\sim\sqrt{\Delta+1}$ where $\Delta$ is a laplacian associated to a smooth Riemannian metric on $M_b$. Choosing a parametrix, the complex powers $Q^{-z}$ are defined for any number $z\in \cc$ with $\re(z)\gg 0$ via an appropriate contour integral (\cite{Se})
\be
Q^{-z}= \frac{1}{2\pi \i} \int_{\Gamma} \la^{-z}(\la-Q)^{-1} d\la\ ,\label{com}
\ee
around the positive real axis. Modulo smoothing operators, we can always arrange $Q^{-z}$ to be a \emph{properly supported} pseudodifferential operator. Taking subsequent products with $Q$ yields the complex powers $Q^{-z}$ for any $z\in \cc$. If $P$ is a \emph{compactly supported} pseudodifferential operator of order $m$, the product $P Q^{-z}$ is a trace-class operator provided that $\re(z)> m+\dim M_b$. The same is true for a product $PU_hQ^{-z}$, where $h$ is a diffomorphism of $M_b$ and $U_h$ the corresponding linear operator on scalar functions. 

\begin{lemma}\label{lzeta}
Let $O\subset G$ be an isotropic submanifold whose action on $M$ is non-degenerate. For any $h\in O$ let $P_{h} $ and $Q_{h}$ be respectively a pseudodifferential operator and an elliptic positive invertible operator of order one, acting along the fiber $M_{r(h)}$ and depending smoothly on the parameter $h$. Then the zeta-function
\be
z\mapsto \Tr(P_{h}U_{h}Q_{h}^{-z})
\ee
defined for $\re(z)\gg 0$ extends to a meromorphic function with simple poles on the complex plane, whose coefficients depend smoothly on $h$.
\end{lemma}
{\it Proof:} The pullback of the submersion $M\to B$ with respect to the map $O\to B$, $h \mapsto r(h)$ is a submersion with base $O$ and fiber $M_{r(h)}$ over any point $h$. Since by hypothesis the submanifold $M_{r(h)}^h$ of fixed points for $h$ varies smoothly with $h$, one can locally choose a system of vertical coordinates $(x;y)$ on this submersion, with the following property: over any point $h$ in a small open set $V\in O$,  $x=(x_1,\ldots,x_r)$ provides a coordinate system on the fixed submanifold $M_{r(h)}^h$, and $y=(y_1,\ldots,y_s)$ provides a normal coordinate system compatible with the diffeomorphism $h$. Then we complete $(x;y)$ into the canonical coordinates $(x,p;y,q)$ on the cotangent bundle $TM_{r(h)}$, such that $(x,p)$ are the canonical coordinates on $TM_{r(h)}^h$ for each $h\in V$. The trace $\Tr(P_{h}U_{h}Q_{h}^{-z})$ is then obtained as the integral, over the manifold $M_{r(h)}$, of a density expressed in the local coordinate system by
$$
\rho^z_{h}(x;y) = \left(\int\!\!\!\int \si^z_{h}(x,p;y,q) \, e^{\i\langle p, x-h^* x\rangle + \i\langle q, y-h^* y \rangle} \, \frac{d^rp d^sq}{(2\pi)^{r+s}}\right) d^syd^rx\ ,
$$
where $\si^z_{h}(x,p;y,q)$ denotes the complete symbol of the pseudodifferential operator $Q_{h}^{-z}P_{h}$, of order $|P|-z$ if $|P|$ is the order of $P$. For notational simplicity we shall drop the subscript $h$ and keep in mind that all symbols depend smoothly on $h$. Note that the symbol $\si^z(x,p;y,q)$ is a holomorphic function of $z$. The above integral converges for $\re(z)\gg 0$ and we want to show that $\rho^z(x;y)$ can be extended to a distribution in the variables $(x,y)$ with values in meromorphic functions of $z\in \cc$. Hence, using local coordinate charts and a partition of unity we get the desired meromorphic extension of the trace. First we perform a change of variables $(x;y)\mapsto (x;u)$ near the submanifold of fixed points, with $u=y-h^*y$. This is allowed because the matrix of partial derivatives $\d u/\d y = 1 - \d(h^*y)/\d y$ is non-singular by hypothesis. Hence the density becomes
$$
\rho^z(x;u) = \left(\int \left(\int \frac{\si^z(x,p;u,q)}{|\det(1-h')|} \, e^{\i\langle p, x-h^* x\rangle} \, \frac{d^rp}{(2\pi)^r}\right) e^{\i\langle q, u \rangle} \, \frac{d^sq}{(2\pi)^s}\right) d^sud^rx
$$
where the matrix $h'=\d(h^*y)/\d y$ is a function of $(x;u)$. The next step is to Taylor expand the symbol $\si^z(x,p;u,q)$ with respect to $q$, up to a certain order $n$. This involves the sequence of holomorphic symbols $\d^k\si^z(x,p;u,q)/\d q^k$ of order $|P|-k-z$:
$$
\si^z(x,p;u,q) = \sum_{k=0}^n \frac{q^k}{k!}\frac{\d^k\si^z}{\d q^k}(x,p;u,0) + \frac{q^{n+1}}{n!}\int_0^1 (1-t)^n \frac{\d^{n+1}\si^z}{\d q^{n+1}} (x,p;u,tq)\, dt\ .
$$
Note that the remainder $R_n^z(x,p;u,q)=\int_0^1 (1-t)^n \frac{\d^{n+1}\si^z}{\d q^{n+1}}(x,p;u,tq)\, dt$ is a not a symbol of order $|P|-n-1-z$ because of integration near $t=0$. Plugging the Taylor expansion in the above expression for $\rho^z(x;u)$, one is left with the terms
\beq
\lefteqn{\frac{1}{k!}\left(\int \left(\int \frac{\d^k\si^z}{\d q^k}(x,p;u,0) \, \frac{e^{\i\langle p, x-h^* x\rangle}}{|\det(1-h')|} \, \frac{d^rp}{(2\pi)^r}\right) q^k e^{\i\langle q, u \rangle} \, \frac{d^sq}{(2\pi)^s}\right) d^sud^rx}\non\\
&&\qquad\qquad = \frac{1}{k!}\left(\int \frac{\d^k\si^z}{\d q^k}(x,p;u,0) \, \frac{e^{\i\langle p, x-h^* x\rangle}}{|\det(1-h')|} \, \frac{d^rp}{(2\pi)^r}\right) \, \frac{\d^k \delta^s(u)}{(\i\d u)^k} d^sud^rx \non
\eeq
where we performed the integral over $q$, and $\delta^s(u)$ is the Dirac mass localized at $(u_1,\ldots,u_s)=(0,\ldots,0)$ which corresponds to the submanifold of fixed points. Hence the $k$-th derivative $\d^k\delta^s(u)/\d u^k$ is a distribution of order $k$ supported by this submanifold. It follows that we only need to know the Taylor expansion around $u=0$ of the integral
$$
I_k=\frac{1}{k!}\int \frac{\d^k\si^z}{\d q^k}(x,p;u,0) \, \frac{e^{\i\langle p, x-h^* x\rangle}}{|\det(1-h')|} \, \frac{d^rp}{(2\pi)^r}
$$
up to order $k$ in the variable $u$, because the higher orders are killed in the product with the $\delta$-distribution. Then the crucial fact is that by the non-degeneracy hypothesis, $x-h^*x$ is of order $u^2$, hence the Taylor expansion of the oscillatory term $\exp(\i\langle p, x-h^* x\rangle)$ yields a \emph{polynomial} in $p$. Therefore $I_k$ reduces to the integral of a classical symbol in the variables $(x,p)$, holomorphic in $z$. By a well-known result it extends to a meromorphic function of $z$ with only simple poles \cite{Wo}. Finally we still have to look at the remainder term
\beq
\lefteqn{\int\!\!\!\int R^z_n(x,p;u,q) \, \frac{e^{\i\langle p, x-h^* x\rangle}}{|\det(1-h')|} \,  q^{n+1} e^{\i\langle q, u \rangle} \, \frac{d^rp}{(2\pi)^r}\frac{d^sq}{(2\pi)^s}}\non\\
&& = \int\!\!\!\int R^z_n(x,p;u,q) \, \frac{e^{\i\langle p, x-h^* x\rangle}}{|\det(1-h')|} \, \frac{\d^{n+1}e^{\i\langle q, u \rangle}}{(\i\d u)^{n+1}}\, \frac{d^rp}{(2\pi)^r}\frac{d^sq}{(2\pi)^s} \non
\eeq
Here we cannot simply perform the integral over $q$ because $R^z_n$ depends on $q$. By the way, this integral will not yield a distribution localized at $u=0$. Instead we shall move the derivatives $\d/\d u$ and rewrite the integral as a sum of terms like
$$
\left(\i\frac{\d}{\d u}\right)^{j}\left(\int\!\!\!\int \left(\i\frac{\d}{\d u}\right)^{m}\left(\frac{R^z_n(x,p;u,q)}{|\det(1-h')|}\right) \frac{\d^{k}e^{\i\langle p, x-h^* x\rangle}}{(-\i\d u)^{k}}  \, e^{\i\langle q, u \rangle}\, \frac{d^rp}{(2\pi)^r}\frac{d^sq}{(2\pi)^s}\right)
$$
with $j+m+k=n+1$. We have $\d^{k}e^{\i\langle p, x-h^* x\rangle}/\d u^{k}=f(x,p;u)e^{\i\langle p, x-h^* x\rangle}$, where the function $f(x,p;u)$ is a polynomial of degree at most $k$ in $p$. Moreover $x-h^*x$ is of order $u^2$ near $u=0$, hence the derivative $\d(x-h^*x)/\d u$ is of order $u$, and the coefficient of $p^l$ in $f(x,p;u)$ is of order $u^{2l-k}$, with $2l-k$ non-negative. But a power of $u$ amounts to a derivative $\d/\d q$ against $e^{\i\langle q, u \rangle}$. Thus we may replace $f$ by a sum of differential operators $p^l (\frac{\d}{\d q})^{2l-k}$ with coefficients smooth functions of $(x;u)$. Each operator amounts to raise the order of $R^z_n(x,p;u,q)$ by $l-(2l-k)=k-l$, which is $\leq k/2$ and hence $\leq (n+1)/2$. Explicitly 
$$
p^l \left(\frac{\d}{\d q}\right)^{2l-k}R^z_n(x,p;u,q)=\int_0^1 (1-t)^nt^{2l-k}\,  p^l\frac{\d^{2l-k+n+1}\si^z}{\d q^{2l-k+n+1}}(x,p;u,tq)\, dt
$$
where $p^l\d^{2l-k+n+1}\si^z/\d q^{2l-k+n+1}$ is a symbol of order $|P|+k-l-n-1-z\leq |P|-\frac{n+1}{2}-z$. Finally, one is left with integrals of the form
$$
J=\left(\i\frac{\d}{\d u}\right)^{j}\left(\int\!\!\!\int\!\!\!\int_0^1 S^z(x,p;u,tq) e^{\i\langle p, x-h^* x\rangle + \i\langle q, u \rangle}L(t)\, dt \frac{d^rp}{(2\pi)^r}\frac{d^sq}{(2\pi)^s}\right)
$$
where $L(t)$ is a polynomial in $t$, and $S^z(x,p;u,q)$ is a holomorphic symbol of order at most $|P|-\frac{n+1}{2}-z$. Thus $S^z$ is dominated by the symbol $(p^2+q^2+1)^{-Z/2}$ for $Z=-|P|+\frac{n+1}{2}+z$. Moreover the integral 
$$
\int\!\!\!\int\!\!\!\int_0^1 (p^2+q^2+1)^{-Z/2} \, dt \frac{d^rp}{(2\pi)^r}\frac{d^sq}{(2\pi)^s}
$$
converges to a holomorphic function of $Z$ provided $\re(Z)$ is large enough. Hence the integral $J$ converges is a holomorphic function of $z$, provided $n$ is chosen sufficiently large. We conclude that $\Tr(PU_hQ^{-z})$ extends to a meromorphic function with simple poles. \\
It remains to show the smoothness with respect to the parameter $h$. In fact it is clear from the integral expression of the density $\rho^z_{h}$ and in all the subsequent calculations, that derivating with respect to $h$ simply amounts to replace the holomorphic symbol $\si^z_{h}$ of order $|P|-z$ by a new holomorphic symbol of order $|P|-z+1$. One concludes that the meromorphic function $\Tr(P_hU_hQ_h^{-z})$ is infinitely differentiable with respect to the parameter $h$. \cqfd\\

We now compute the residue at $z=0$ of a zeta-function of type $\Tr(PU_hQ^{-z})$ with $P$ a compactly supported pseudodifferential operator on the manifold $M_b$, $h$ a diffeomorphism of $M_b$ and $Q$ an elliptic positive operator of order one. Not surprisingly, the residue is given by an explicit local formula involving the complete symbol $\si_P$ of $P$. In particular when $U_h$ is the identity, one recovers the well-known Wodzicki residue \cite{Wo}, which can be written as an integral, over the cosphere bundle, of a certain homogeneous component of $\si_P$. In the general situation the residue is localized at the set of fixed points for $U_h$. 

\begin{proposition}\label{ploc}
Let $h$ be a diffeomorphism of the manifold $M_b$. Assume that the set of fixed points of $h$ is a non-degenerate smooth submanifold $M_b^h\subset M_b$ of dimension $r$. Choose a local coordinate system $(x;y)$ adapted to $M_b^h$, and complete it into a canonical coordinate system $(x,p;y,q)$ on the cotangent bundle $T^*M_b$, such that $(x,p)$ is a canonical coordinate system of $T^*M_b^h$. Then for any pseudodifferential operator $P\in \CL_c^k(M_b)$ and any elliptic strictly positive invertible operator $Q\in \CL^1(M_b)$, one has the localization formula
\be
\res \Tr(PU_h Q^{-z}) = \int_{S^*M_b^h} \left[ e^{\i\langle \frac{\d}{\d q}, (1-h')^{-1}\frac{\d}{\d y}\rangle} \cdot  \frac{ \si_P\, e^{\i\langle p, x-h^*x\rangle} }{|\det(1-h')|}  \right]_{-r} \, \frac{\eta (d\eta)^{r-1}}{(2\pi)^r} \label{floc}
\ee
where $S^*M_b^h$ is the cosphere bundle of the fixed submanifold, $\eta=\langle p , dx\rangle$ is the canonical one-form on the cotangent bundle $T^*M_b^h$, $\si_P=\si_P(x,p;y,q)$ is the complete symbol of $P$, $[\ ]_{-r}$ is the order $-r$ component of a symbol in the variables $(x,p; 0,0)$, and $h'$ is the matrix of partial derivatives $\d (h^*y)/\d y$.  
\end{proposition}
{\it Proof:} Set $u=y-h^*y$. In the proof of Lemma \ref{lzeta} we established that $\Tr(PU_hQ^{-z})$ is the integral of a density on $M_b$ given in the local coordinate chart $(x;u)$ by an expansion
\beq
\rho^z(x;u) &=& \sum_{k=0}^n \frac{1}{k!}\left(\int \frac{\d^k\si^z}{\d q^k}(x,p;u,0) \, \frac{e^{\i\langle p, x-h^* x\rangle}}{|\det(1-h')|} \, \frac{d^rp}{(2\pi)^r}\right) \, \frac{\d^k \delta^s(u)}{(\i\d u)^k} d^sud^rx \non\\
&& + H_n(z)\non
\eeq
where the remainder $H_n(z)$ is holomorphic in $z$ provided $n$ is sufficiently large, and $\si^z$ is the symbol of $Q^{-z}P$ of order $|P|-z$. Hence integrating over $(x;u)$ and taking the residue at $z=0$ will only retain the finite sum over $k$:
$$
\res \sum_{k=0}^n \frac{1}{k!}\int\!\!\!\int\!\!\!\int \frac{\d^k\si^z}{\d q^k}(x,p;u,0) \, \frac{e^{\i\langle p, x-h^* x\rangle}}{|\det(1-h')|} \, \frac{\d^k \delta^s(u)}{(\i\d u)^k} \, \frac{d^rp}{(2\pi)^r}d^sud^rx 
$$
Integrating by parts with respect to the variable $u$ one sees that the Dirac measure $\delta$ localizes the residue at the submanifold of fixed points $u=0$:
$$
\res \int\!\!\!\int \left(\sum_{k=0}^n \frac{1}{k!} \left\langle \frac{\d^k}{\d q^k},\frac{\d^k}{(-\i\d u)^k}\right\rangle \right)\cdot \left. \frac{\si^ze^{\i\langle p, x-h^* x\rangle}}{|\det(1-h')|}\right|_{\substack{u=0 \\ q=0}} \, \frac{d^rp}{(2\pi)^r}d^rx 
$$
Moreover we know that $x-h^*x$ is of order $u^2$ near $u=0$, and the Taylor expansion of $\d^k\si^z/\d q^k(x,p;u,0) \, e^{\i\langle p, x-h^* x\rangle}$ up to order $k$ in the variable $u$ is a symbol with respect to the variables $(x,p)$, of order at most $|P|-k/2-z$. Hence for high values of $k$ the integral over $(x,p)$ converges to a holomorphic function of $z$ near $z=0$ and the residues vanish. Consequently we may replace the finite sum $\sum_{k=0}^n \frac{1}{k!} \langle \frac{\d^k}{\d q^k},\frac{\d^k}{(-\i\d u)^k}\rangle$ by the exponential $\exp(\i\langle \frac{\d}{\d q},\frac{\d}{\d u}\rangle)$ viewed as a formal power series. Moreover, the symbol $\si^z$ of the operator $Q^{-z}P$ has an asymptotic expansion of the form
$$
\si^z\sim (\si_Q)^{-z}\si_P + z \si'
$$
where $\si_Q$ is the symbol of $Q$ and the product is simply the product of functions in the canonical variables. The remainder $z\si'$ will disappear under the residue because the integral of a symbol with respect to $(x,p)$ is meromorphic with only simple poles. Hence we can replace $\si^z$ by the product $(\si_Q)^{-z}\si_P$. In the same way we can move the function $(\si_Q)^{-z}$ behind the differential operators $\langle \frac{\d^k}{\d q^k},\frac{\d^k}{\d u^k}\rangle$ since the latter can only extract symbols proportional to $z$ when applied to $(\si_Q)^{-z}$. Thus the above integral may be rewritten
$$
\res \int\!\!\!\int (\si_Q)^{-z}\, e^{\i\langle \frac{\d}{\d q}, \frac{\d}{\d u}\rangle} \cdot \left. \frac{\si_Pe^{\i\langle p, x-h^* x\rangle}}{|\det(1-h')|}\right|_{\substack{u=0 \\ q=0}} \, \frac{d^rp}{(2\pi)^r}d^rx 
$$
One recognizes the Wodzicki residue applied to pseudodifferential operator with symbol $e^{\i\langle \frac{\d}{\d q}, \frac{\d}{\d u}\rangle} \cdot\si_Pe^{\i\langle p, x-h^* x\rangle}/|\det(1-h')|_{u=0,q=0}$ on the cotangent bundle of the fixed submanifold $M_b^h$. It is expressed as the integral of the order $-r$ component of the symbol over the cosphere bundle, whence formula (\ref{floc}).  \cqfd\\

Following \cite{P8}, the zeta-renormalized trace of an operator of the form $P U_h$ can be defined as the \emph{finite part} at $z=0$ of the zeta-function $\Tr(PU_hQ^{-z})$, that is, the term of degree zero in the Laurent expansion of the zeta-function at $z=0$. Choose a generalized connection $h$ on the submersion $\pi: M\to B$ and a smooth section $Q\in \cinf(B,\CL^1(M))$ of elliptic, positive and invertible pseudodifferential operators as above. Thus at any point $b\in B$ one has an elliptic positive invertible operator $Q_b$ acting on the manifold $M_b=\pi^{-1}(b)$. Let $O\subset G$ be an $\Ad$-invariant isotropic submanifold whose action on $M$ is non-degenerate. Then by inserting the complex power $Q^{-z}$ inside the trace map (\ref{toto}), we associate to any tensor $e_1\otimes \ldots \otimes e_n\in T\Ec$ the function $\Tr^{h,Q}_{e_1, \ldots , e_n}(z)\in C^k(V)$ defined over an appropriate neighborhood $V$ of the submanifold $O^{(n)}\subset G^n$:
\beq
\lefteqn{\Tr^{h,Q}_{e_1, \ldots , e_n}(z)(g_1,\ldots, g_n) = }  \label{trzeta0} \\ 
&& \qquad \Tr\big(e_1(g_1)\, U_{g_1}\, {h(s(g_1),r(g_2))}  \ldots e_n(g_n)\, U_{g_n}\, {h(s(g_n),r(g_1))}\, Q^{-z}_{r(g_1)}\big) \non
\eeq
The regularity order $k$ can be as large as wanted, provided that $\re(z)$ is large enough. Hence by Lemma \ref{lzeta}, the function $\Tr^{h,Q}_{e_1, \ldots , e_n}(z)$ \emph{projected to the localization space $\cinf(V)_{[O]}$ of jets to all order at $O$} extends to a meromorphic function of $z$ with at most simple poles. This implies that the evaluation of $\Tr^{h,Q}_{e_1, \ldots , e_n}(z)$ on a distribution $\varphi\in C^{-\infty}(V)$ with support localized at $O$ and of bounded singularity order, yields a meromorphic function of $z$. We also introduce the notation
\beq
\lefteqn{\Tr^{h,Q}_{e_1, \ldots , e_n}(z)(g_1,\ldots, g_{n-1}|g_n) = }  \label{trzeta1} \\ 
&& \qquad \Tr\big(e_1(g_1)\, U_{g_1}\, {h(s(g_1),r(g_2))}  \ldots  Q^{-z}_{r(g_n)}\,  e_n(g_n)\, U_{g_n}\, {h(s(g_n),r(g_1))}\big) \non
\eeq
Since we regard $\Tr^{h,Q}_{e_1, \ldots , e_n}(z)$ as the jets of a function at the submanifold $O^{(n)}$, we can take any coefficient of its Laurent expansion at $z=0$ \emph{before} evaluating it on a localized distribution.

\begin{definition}
Let $\Cc = \cinfc(B)\cp G$ and $\Ec = \cinfc(B,\CL^0_c(M))\cp G$. Choose any generalized connection $h$ on the submersion $\pi:M\to B$ and any elliptic section $Q$ as above. Then for any cocycle $\varphi\in \hom(X(\Th\Cc_\top)_{[O]},\cc)$ localized at $O$, the zeta-renormalized cochain $(\tau_\varphi)_R \in \hom(X(\Th\Ec),\cc)$ is defined on $T\Ec$ and $\Om^1T\Ec_\nat$ by
\beq
&&\hspace{-1cm} (\tau_\varphi)_R(e_1\otimes\ldots\otimes e_n) = \int_{G^n} \varphi^+_n(g_1,\ldots,g_n)\, \Pf_{z=0} \Tr^{h,Q}_{e_1,\ldots,e_n}(z)(g_1,\ldots,g_n)\ , \non \\
&&\hspace{-1cm}(\tau_\varphi)_R(\nat(e_1\otimes\ldots\otimes e_{n-1}\dd e_n)) = \label{zeta} \\
&&\qquad\qquad \int_{G^n} \varphi^-_n(g_1,\ldots,g_{n-1}|g_n)\, \Pf_{z=0}\Tr^{h,Q}_{e_1,\ldots,e_n}(z)(g_1,\ldots,g_{n-1}|g_n) \ .\non
\eeq
where $\Pf_{z=0}$ denotes the finite part at $z=0$ of the corresponding zeta-functions.
\end{definition}
Observe that $(\tau_\varphi)_R$ is well-defined on the pro-complex $X(\Th\Ec)$ because evaluation on the distributions $\varphi^\pm$ kills the high powers of the ideal $J\Ec$. The zeta-function also allows to define a \emph{residue morphism}
\be 
\Res\ :\ X(\Tt\Ec) \to X(\Th\Cc_\top)_{[O]} \label{resmap}
\ee
by selecting the poles of (\ref{trzeta0}) and (\ref{trzeta1}) at $z=0$ instead of the finite part. The restriction of (\ref{resmap}) to the even subspace of the $X$-complex is a linear map $\Tt\Ec\to (\Th\Cc_\top)_{[O]}$, sending any $n$-tensor $e_1\otimes\ldots \otimes e_n \in T\Ec$ to the jets of a function of the variables $(g_1,\ldots,g_n)$ at the localization submanifold $O^{(n)}\subset G^n$: 
$$
\big(\Res(e_1\otimes\ldots \otimes e_n)\big)(g_1,\ldots,g_n) = \res \Tr^{h,Q}_{e_1, \ldots , e_n}(g_1,\ldots, g_n)(z)\ .
$$
Notice that the ideal $J\Ec\subset T\Ec$ is sent to $(J\Cc_\top)_{[O]}$. Moreover, for any fixed $k$ the $k$-jet of $\Res(e_1\otimes\ldots \otimes e_n)$ vanishes whenever $e_1\otimes\ldots \otimes e_n$ belongs to a sufficiently high power of the ideal $\Rc=\ker(T\Ec\to T\Ac)$. This is due to the fact that the zeta-function (\ref{trzeta0}) has no pole at $z=0$ in this case. Hence the residue morphism indeed extends to a well-defined linear map $\Tt\Ec\to (\Th\Cc_\top)_{[O]}$. In the odd case (\ref{resmap}) is a linear map $\Om^1\Tt\Ec_\nat\to (\Om^1\Th\Cc_\top)_{[O]}$ defined in an analogous way:
$$
\big(\Res(\nat(e_1\otimes\ldots \otimes e_{n-1}\dd e_n))\big)(g_1,\ldots,g_{n-1}|g_n) = \res \Tr^{h,Q}_{e_1, \ldots , e_n}(g_1,\ldots, g_{n-1}|g_n)(z)\ .
$$
The crucial point is that (\ref{resmap}) is a morphism of $X$-complexes. This is an easily consequence of the fact that the zeta-functions have only simple poles: the residues at $z=0$ do not depend on the actual place of $Q^{-z}$ in formulas (\ref{trzeta0}) and(\ref{trzeta1}). This property would fail in the presence of double poles. Following \cite{P8} we now introduce the logarithm 
\be 
\ln Q = - \frac{d}{dz}Q^{-z}|_{z=0}\ .
\ee
The latter is no longer a section of the classical pseudodifferential operators $\CL(M)$, but belongs to the larger class of \emph{log-polyhomogeneous} pseudodifferential operators $\CL(M)_{\log}$. However, the difference $\ln Q-\ln Q'$ of two such logarithms \emph{is} a section of $\CL(M)$, as well as the commutator $[\ln Q, P]$ with any section $P\in \cinfc(B,\CL_c(M))$. We shall enlarge the algebra $\Ec$ by adding the log-polyhomogeneous operators. Define
$$
\Ec_{\log} = \cinfc(B,\CL^0_c(M)_{\log})\cp G\ .
$$
Thus the elements of $\Ec_{\log}$ are products of logarithms by elements of $\Ec$. In particular the section $\ln Q$ can be viewed as a left multiplier of $\Ec_{\log}$, as follows: $(\ln Q\cdot e)(g) = \ln Q_{r(g)} e(g)$ for all $e\in \Ec$ and $g\in G$. We can mimic the construction of the complex $X(\Tt\Ec)$, replacing everywhere classical pseudodifferential operators by log-polyhomogeneous ones. This leads to a complex $X(\Tt\Ec)_{\log}$ and its subcomplex
\be 
X(\Tt\Ec)_{\log}^1 \subset X(\Tt\Ec)_{\log}\ ,
\ee
where the superscript $^1$ means that we retain only the tensors having logarithmic degree at most $1$. Thus, the elements of $(T\Ec)_{\log}^1$ are of the form $e_1\otimes \ldots \otimes e_n$ or $e_1\otimes \ldots \otimes \ln Q \cdot e_i\otimes \ldots \otimes e_n$, where all $e_j$'s are in $\Ec$. Similarly in odd degree, the elements of $(\Om^1T\Ec_\nat)_{\log}^1$ are of the form $\nat( e_1\otimes \ldots \otimes e_{n-1}\dd e_n)$ , $\nat(e_1\otimes \ldots \ln Q \cdot e_i \ldots \otimes e_{n-1}\dd e_n)$ or $\nat( e_1\otimes \ldots \otimes e_{n-1}\dd(\ln Q\cdot e_n))$. Since the difference of logarithms $\ln Q_{r(g)} - U_g\ln Q_{s(g)}U_g^{-1}$ is always a classical pseudodifferential operator on the manifold $M_{r(g)}$, one sees that the residue map (\ref{resmap}) can be extended to a chain map
$$
X(\Tt\Ec)^1_{\log}\cap\dom(\Res) \to X(\Th\Cc_\top)_{[O]}
$$
where the domain $\dom(\Res)$ is the linear span of differences of chains where only the place of $\ln Q$ changes. For example in even degree, the chains in $(T\Ec)^1_{\log}$ are linearly generated by differences
$$
e_1\otimes\ldots \ln Q\cdot e_i\otimes \ldots e_j\otimes  \ldots \otimes e_n - e_1\otimes\ldots e_i\otimes \ldots \ln Q\cdot e_j\otimes \ldots \otimes e_n
$$
For notational convenience we introduce the convention that a right multiplication of a factor $e_i$ by $\ln Q$ amounts to the left multiplication of the following factor $e_{i+1}$ by $\ln Q$ in a tensor product. In particular
\beq
\lefteqn{e_1\otimes\ldots \otimes[\ln Q, e_i]\otimes e_{i+1}\otimes  \ldots \otimes e_n := }\non\\
&&  e_1\otimes\ldots \otimes \ln Q \cdot e_i \otimes e_{i+1}\otimes  \ldots \otimes e_n - e_1\otimes\ldots \otimes e_i \otimes \ln Q\cdot e_{i+1}\otimes  \ldots \otimes e_n \non\\
&&\hspace{-0.7cm} e_1\otimes\ldots \otimes e_n\cdot \ln Q := \ln Q\cdot e_1\otimes\ldots \otimes e_n \non
\eeq

\begin{proposition}\label{pres}
Let $\Cc = \cinfc(B)\cp G$ and $\Ec = \cinfc(B,\CL^0_\pi(M))\cp G$. Choose any generalized connection $h$ on the submersion $\pi:M\to B$ and any elliptic section $Q$ as above. Then for any cocycle $\varphi\in \hom(X(\Th\Cc_\top)_{[O]},\cc)$ localized at $O$, the  boundary of the zeta-renormalized cochain $(\tau_\varphi)_R$ is the cocycle $(\tau_\varphi)_R\d \in \hom(X(\Tt\Ec),\cc)$ given by 
\beq
&&\hspace{-1cm} (\tau_\varphi)_R\d(\nat( e_1\otimes \ldots\otimes e_{n-1}\dd e_n)) = \varphi\circ\Res(e_1\otimes \ldots\otimes e_{n-1} \otimes  [\ln Q, e_n]) \label{bound} \\
&&\hspace{-1cm} (\tau_\varphi)_R\d(e_1\otimes \ldots\otimes e_n) =  \sum_{1\leq i<j\leq n} \varphi\circ\Res( e_1 \otimes \ldots [\ln Q, e_i] \ldots \dd e_j \ldots \otimes e_n) \non
\eeq
\end{proposition}
{\it Proof:} By definition one has $\d (\nat( e_1\otimes \ldots\otimes e_{n-1}\dd e_n)) = e_1\otimes \ldots\otimes e_{n-1} \otimes e_n - e_n\otimes e_1\otimes \ldots\otimes e_{n-1}$. For all $i$ write $h^i_{i+1} = h(s(g_i),r(g_{i+1}))$. Then (\ref{zeta}) gives
\beq
\lefteqn{(\tau_\varphi)_R(e_1\otimes \ldots\otimes e_{n-1}\otimes e_n) =  \int_{G^n} \varphi^+_n(g_1,\ldots,g_n)\times } \non\\
&& \Pf_{z=0}\Tr\big( e_1(g_1) U_{g_1} h^1_2 \ldots e_{n-1}(g_{n-1}) U_{g_{n-1}} h^{n-1}_n e_n(n_n) U_{g_n} h^n_1 Q_{r(g_1)}^{-z}\big) \non 
\eeq
In the same way
\beq
\lefteqn{(\tau_\varphi)_R(e_n\otimes e_1\otimes \ldots\otimes e_{n-1}) = \int_{G^n} \varphi^+_n(g_1,\ldots,g_n)\times  } \non\\
&& \Pf_{z=0}\Tr\big( e_n(g_1) U_{g_1} h^1_2 e_1(g_2) U_{g_2} h^2_3 \ldots  e_{n-1}(g_n) U_{g_n} h^n_1 Q_{r(g_1)}^{-z}\big)\ . \non 
\eeq
Since $\varphi$ is an $X$-complex cocycle, $\varphi^+_n(g_1,\ldots,g_n)$ is invariant under cyclic permutations of $(g_1,\ldots,g_n)$. This and the cyclicity of the operator trace implies
\beq
\lefteqn{(\tau_\varphi)_R(e_n\otimes e_1\otimes \ldots\otimes e_{n-1}) =  \int_{G^n} \varphi^+_n(g_1,\ldots,g_n)\times } \non\\
&& \Pf_{z=0}\Tr\big( e_1(g_1) U_{g_1} h^1_2 \ldots e_{n-1}(g_{n-1}) U_{g_{n-1}} h^{n-1}_n Q_{r(g_n)}^{-z} e_n(g_n) U_{g_n} h^n_1 \big) \non 
\eeq
Thus one can write
\beq
\lefteqn{(\tau_\varphi)_R\d(\nat(e_1\otimes \ldots\otimes e_{n-1} \dd e_n))= \int_{G^n} \varphi^+_n(g_1,\ldots,g_n)\times } \non\\
&& \Pf_{z=0}\Tr\big( e_1(g_1) U_{g_1} h^1_2 \ldots e_{n-1}(g_{n-1}) U_{g_{n-1}} h^{n-1}_n [e_n(g_n) U_{g_n} h^n_1, Q^{-z}]\big) \non 
\eeq
with the ``commutator'' 
\beq
\lefteqn{[e_n(g_n) U_{g_n} h^n_1, Q^{-z}]= e_n(g_n) U_{g_n} h^n_1 Q_{r(g_1)}^{-z} - Q_{r(g_n)}^{-z} e_n(g_n) U_{g_n} h^n_1 } \non\\
&=& e_n(g_n) \big(U_{g_n} h^n_1 Q_{r(g_1)}^{-z} (U_{g_n} h^n_1)^{-1}- Q_{r(g_n)}^{-z} \big) U_{g_n} h^n_1 - [Q_{r(g_n)}^{-z}, e_n(g_n)] U_{g_n} h^n_1 \non
\eeq
Now observe that $U_{g_n} h^n_1 Q_{r(g_1)}^{-z}(U_{g_n} h^n_1)^{-1}- Q_{r(g_n)}^{-z}$ and $[Q_{r(g_n)}^{-z}, e_n(g_n)]$ are pseudodifferential operators on the manifold $M_{r(g_n)}$. They have an asymptotic expansion in powers of $z$,
\beq
\lefteqn{U_{g_n} h^n_1 Q_{r(g_1)}^{-z} (U_{g_n} h^n_1)^{-1}- Q_{r(g_n)}^{-z}} \non\\
&&\quad \sim -z \big( U_{g_n} h^n_1 \ln(Q_{r(g_1)}) Q_{r(g_1)}^{-z} (U_{g_n} h^n_1)^{-1}- \ln(Q_{r(g_n)}) Q_{r(g_n)}^{-z} \big) + O(z^2)\non\\
&&\hspace{-0.5cm} {[Q_{r(g_n)}^{-z}, e_n(g_n)]} \sim  -z [\ln Q_{r(g_n)},e_n(g_n) ]Q_{r(g_n)}^{-z} + O(z^2) \non
\eeq
up to order $z^2$. Hence with obvious notations
$$
[e_n(g_n) U_{g_n} h^n_1, Q^{-z}] \sim z [\ln Q, e_n(g_n)U_{g_n} h^n_1] Q_{r(g_1)}^{-z} + O(z^2)\ .
$$
Because the zeta-functions have only simple poles, $\Pf_{z=0}(z\Tr(\ldots Q^{-z}))$ is the residue of $\Tr(\ldots Q^{-z})$ and the terms of order $z^2$ are killed. Finally
\beq
\lefteqn{(\tau_\varphi)_R\d(\nat(e_1\otimes \ldots\otimes e_{n-1}\dd e_n)) =  \int_{G^n} \varphi^+_n(g_1,\ldots,g_n)\times } \non\\
&& \res\Tr\big( e_1(g_1) U_{g_1} h^1_2 \ldots e_{n-1}(g_{n-1}) U_{g_{n-1}} h^{n-1}_n [\ln Q, e_n(g_n) U_{g_n} h^n_1] Q_{r(g_1)}^{-z}\big) \non 
\eeq
which is precisely (\ref{bound}). One proceeds similarly with the second formula. \cqfd\\

Collecting the preceding results one gets the following refinement of Theorem \ref{ttop} which computes the excision map by means of a residue formula.

\begin{theorem}\label{tres}
Let $G\rightrightarrows B$ be a Lie groupoid and let $O$ be an $\Ad$-invariant isotropic submanifold of $G$. Let $\pi:M\to B$ be a $G$-equivariant surjective submersion and assume the action of $O$ on $M$ non-degenerate. Then one has a commutative diagram 
\be 
\vcenter{\xymatrix{
HP^\bullet(\cinfc(B,\CL_{c}^{-1}(M))\cp G) \ar[r]^{\quad E^*} & HP^{\bullet+1}(\cinfc(S^*_\pi M)\rtimes G)  \\
HP_\top^{\bullet}(\cinfc(B)\cp G)_{[O]} \ar[u]^{\tau_*} \ar[r]^{\pi^!_G\qquad} & HP_\top^{\bullet+1}(\cinfc(S^*_\pi M)\rtimes G)_{[\pi^* O]} \ar[u] }}
\ee
where the isotropic submanifold $\pi^*O \subset S^*_\pi M\rtimes G$ is the pullback of $O$ by the submersion $S^*_\pi M\to B$.\\
Let $\Ac = \cinfc(S^*_\pi M)\cp G$, $\Ec = \cinfc(B,\CL_c^0(M))\cp G$, and choose any continuous linear splitting $\si:\Ac\to\Ec$ of the projection homomorphism. Then the image of an even class $[\varphi]\in HP_\top^{0}(\cinfc(B)\cp G)_{[O]}$ is represented by the odd cyclic cocycle $\pi^!_G(\varphi)\in \hom(\Om^1\Th\Ac_\nat,\cc)$ over the algebra $\Ac=\cinfc(S^*_\pi M\cp G)$, given by the residue
$$
\pi^!_G(\varphi)(\nat (a_1\otimes \ldots\otimes a_{n-1} \dd a_n)) =  \varphi\circ\Res\big( \si(a_1)\otimes \ldots \otimes\si(a_{n-1})\otimes [\ln Q, \si(a_n)] \big)
$$
for all $\nat(a_1\otimes \ldots\otimes a_{n-1}\dd a_n) \in \Om^1 T\Ac_\nat$. In a similar way, the image of an odd class $[\varphi]\in HP_\top^{1}(\cinfc(B)\cp G)_{[O]}$ is represented by the cyclic cocycle of even degree $\pi^!_G(\varphi)\in \hom(\Th\Ac,\cc)$ given by the residue 
$$
\pi^!_G(\varphi)(a_1\otimes \ldots\otimes a_n) = \!\!\!  \sum_{1\leq i<j\leq n} \!\!\! \varphi\circ\Res\big( \si(a_1)\otimes \ldots [\ln Q, \si(a_i)]\ldots \dd\si(a_j)\ldots \otimes\si(a_n)\big) 
$$
for all $a_1\otimes \ldots\otimes a_n \in T\Ac$. \cqfd
\end{theorem}

Corollary \ref{cind} and Theorem \ref{tres} allow to compute the pairing $\langle [\tau_\varphi],\Ind_E([u])\rangle$ for any elliptic symbol $u\in GL_{\infty}(\Ac)$. For clarity we suppose $u\in GL_1(\Ac)$ but the general case follows easily. Hence let $\si(u)\in \Ec^+$ be the image of $u$ under the linear splitting $\si:\Ac\to\Ec$ (extended to unitalized algebras so that $\si(1)=1$). Take the invertible lifting $\uh=u$ in $\Th\Ac^+$. The unitalized homomorphism $\si_*:\Th\Ac^+\to\Tt\Ec^+$ carries $\uh$ to an invertible element $\si_*(\uh)$, whose inverse is given by the convergent series
$$
\si_*(\uh^{-1})=\sum_{n=0}^{\infty} \si(u^{-1})\otimes(1-\si(u)\otimes\si(u^{-1}))^{\otimes n}\ \in GL_1(\Tt\Ec)\ .
$$

\begin{corollary}\label{ctr}
Under the hypotheses of Theorem \ref{tres}, let $[u]\in K_1(\Ac)$ be an elliptic symbol class element represented by an invertible matrix $u\in M_{\infty}(\Ac)^+$. Let $[\varphi]\in HP_\top^{0}(\cinfc(B)\cp G)_{[O]}$ be a cyclic cohomology class localized at $O$. Then
\be
\langle [\tau_\varphi] , [u] \rangle =   \varphi\circ\Res\#\tr\big( \si_*(\uh^{-1}) [\ln Q, \si_*(\uh)] \big)
\ee
where $Q\in \cinf(B,\CL^1(M))$ is any section of elliptic positive invertible pseudodifferential operators of order one, and $\uh\in M_{\infty}(\Th\Ac)^+$ is the above invertible lifting of $u$. \cqfd
\end{corollary}

\end{document}